\documentstyle[amssymb,12pt]{amsart}

\newtheorem{proposition}{Proposition}[section]
\newtheorem{lemma}{Lemma}[section]
\newtheorem{theorem}{Theorem}[section]
\newtheorem{corollary}{Corollary}[section]

\newcounter{remark}[section]
\renewcommand{\theremark}{\arabic{section}.\arabic{remark}}
\newenvironment{remark}%
{\refstepcounter{remark}\trivlist \item[\hskip
    \labelsep {\bf Remark \theremark}]}%
{\endtrivlist}

\begin{document}
\newcommand{\nc}{\newcommand}
\nc{\pa}{\partial}
\nc{\cA}{{\cal A}}\nc{\cB}{{\cal B}}\nc{\cC}{{\cal C}}
\nc{\cE}{{\cal E}}\nc{\cG}{{\cal G}}\nc{\cH}{{\cal H}}\nc{\cZ}{{\cal
Z}}\nc{\cW}{{\cal W}} 
\nc{\cX}{{\cal X}}
\nc{\cY}{{\cal Y}}
\nc{\cR}{{\cal R}}\nc{\Ad}{\on{Ad}}
\nc{\ad}{\on{ad}}\nc{\Der}{\on{Der}}\nc{\End}{\on{End}}
\nc{\Imm}{\on{Im}}
\nc{\de}{\delta}\nc{\si}{\sigma}\nc{\on}{\operatorname}
\nc{\al}{\alpha}
\nc{\CC}{{\Bbb C}}\nc{\ZZ}{{\Bbb Z}}
\nc{\la}{{\lambda}}
\nc{\wt}{\widetilde}
\nc{\G}{\wt{\g}}
\nc{\g}{{\frak g}}
\nc{\h}{{\frak h}}
\nc{\n}{{\frak n}}
\nc{\ep}{\epsilon}

\title[Geometric interpretation of Poisson structure in Toda theories]
{Geometric interpretation of the Poisson structure in affine Toda field
theories}\thanks{The research of the second author was partially supported
by grants from the Packard Foundation, NSF and the Sloan Foundation}

\author{Benjamin Enriquez}
\address{Centre de Math\'ematiques, URA 169 du CNRS, Ecole Polytechnique,
91128 Palai-seau, France}

\author{Edward Frenkel}
\address{Department of Mathematics, Harvard University,
Cambridge, MA 02138, USA}

\date{June 1996}
\maketitle

\begin{abstract} 
We express the Poisson brackets of local fields of the affine Toda field
theories in terms of the Drinfeld-Sokolov dressing operator. For this, we
introduce a larger space of fields, containing ``half screening charges''
and ``half integrals of motions''. In addition to local terms, the Poisson
brackets contain nonlocal terms related to trigonometric $r$-matrices.
\end{abstract}

\section*{Introduction}

Since the work of Zakharov and Shabat \cite{ZS}, the dressing
techniques have played an important role in the theory of classical
integrable systems. These techniques have been developed by Drinfeld
and Sokolov in \cite{DS} in the framework of affine Toda field
theories. Later, Feigin and one of us proposed (\cite{FF-CIME},
\cite{FF-Inv}) another approach to these theories; this approach was
shown (\cite{E-TMP}, \cite{EFr}) to be equivalent to that of
\cite{DS}. In those works, the space of local fields of the Toda
theory (equivalently, the mKdV hierarchy) associated to an affine Lie
algebra ${\frak g}$ is described as the ring of functions on the coset
space $N_{+}/A_{+}$ of a unipotent subgroup of the Kac-Moody group $G$
corresponding to ${\frak g}$. The mKdV flows are then identified with
the right action of the principal commutative Lie algebra ${\frak a}$
normalizing $A_{+}$, $N_{+}$ being viewed as an open subset of the
flag manifold of $G$. This leads to a system of variables, in which
the flows become linear and hence can be integrated.

In the works on quantization of the Toda theories, an important role
is played by the vertex operator algebra structure on the space of
local fields. At the classical level, this gives rise to what we call
here a vertex Poisson algebra (VPA) structure on the space of local
fields of a Toda theory. The notion of the VPA structure coincides
with the notion of ``coisson algebra'' (on the disc) introduced by
Beilinson and Drinfeld in \cite{BD}. The goal of this work is to
define this and related structures on the space of fields of a Toda
theory in the Lie group terms using the identification described
above.

To this end, we make use of an idea introduced earlier by Feigin and one of
us in \cite{EFe}, where a similar problem was solved in the setting of the
classical lattice Toda theory associated to $\widehat{\frak{sl}}_2$. In
that work, the space of local fields was extended by the ``half screening
charges'' and ``half integrals of motions''. The screening charges and
integrals of motions (IM's) are sums of the lattice translates of certain
expressions, and their ``halfs'' are just the sums of positive translates
of the same expressions. The full space was then identified with the
quotient $G/H$, where $H$ is the Cartan subgroup of $G$. The smaller spaces
of fields, without half screening charges (resp. without half IM's), can be
obtained by taking the quotient of $G/H$ from the left by a Borel subgroup
$B_{-}$ (resp. from the right by the positive part of the loop group of the
Cartan subgroup).  In this interpretation, the Poisson structure on the
full space is given by the difference $\ell(R)-r(R^{\infty})$, where
$\ell(R)$ stands for the left action of the trigonometric $r$-matrix $R\in
{\frak g}^{\otimes 2}$ and $r(R^{\infty})$ stands for the right action of
its ``infinitely twisted'' version $R^{\infty}\in {\frak g}^{\otimes
2}$. This leads to a description of the Poisson structures of the
quotients.

In this work, we enlarge the space of local fields of the continuous affine
Toda theories in a similar way, by adding continuous analogues of the half
screening charges and IM's. Using results of \cite{EFe,EFr}, we identify
the full space $\bar\pi_{0}$ obtained this way with the space of functions
on $B_{-}\times N_{+}$. We then study the structure of nonlocal VPA on
$\bar\pi_0$. The axioms for this structure are given in Section 1. The main
feature is that the the Poisson bracket $\{u(x),v(y)\}$, where $u,v$ are
elements of the algebra, can be expressed as a linear combination of local
terms of the form $u_{k}(x)\pa_{x}^{k}\delta(x-y)$ with $k\ge 0$, and of
nonlocal terms of the form $a(x)\pa_{x}^{-1}\delta(x-y)b(y)$.

A similar formalism was introduced by Radul \cite{Radul} in the framework
of formal variational calculus (see also \cite{LR}). Natural examples of
nonlocal VPA's are given by the higher Adler-Gelfand-Dickey (AGD)
structures (i.e. structures obtained from the pair of the first two AGD
structures by application of the Magri recursion procedure) introduced in
\cite{EOR,AvM}.

In Sections 2 and 3 we give a geometric description of the nonlocal
VPA structure of $\bar\pi_{0}=\CC[B_{-}\times N_{+}]$. A generic
element $g=(b_{-},n_{+})$ of $B_{-}\times N_{+}$ can be considered as
the product of expansions of the scattering matrix of the Lax operator
from $-\infty$ to $x$ for a small spectral parameter $\lambda$, and
from $x$ to $+\infty$ for a large $\lambda$.

The Poisson bracket on $N_+$ is obtained by a straightforward
extension of the local VPA structure on $N_+/A_+$. On the other hand,
according to \cite{BLZ} and in the spirit of \cite{Faddeev},
\cite{Bab}, we define the Poisson bracket on $B_{-}$ via the
trigonometric $r$-matrix. We show that the trigonometric Poisson
brackets on $B_{-}$ are compatible with the Poisson brackets of local
fields (see Lemma \ref{lemma2.2}). The Poisson brackets of $N_{+}$ and
between $N_+$ and $B_-$ also have nonlocal parts which we determine in
Lemmas \ref{lemma3.4} and \ref{lemma3.5}. To derive the complete
expression for the Poisson brackets (Thm. \ref{thm3.1}), we use the
evolution equation
\begin{equation}    \label{gp}
\pa_{x}g(x)=g(x)p_{-1},
\end{equation}
where $p_{-1}$ is a degree $-1$ element of ${\frak a}$.

After that we obtain another realization of the VPA structures on $N_{+}$
and $N_{+}/A_{+}$ purely in terms of the unipotent group elements (i.e. the
Drinfeld-Sokolov dressing operators) -- see Cor. 3.1 and formula
(\ref{finalPBn-n}). The latter formula could in principle be obtained
directly in the framework of $N_+$, but the simple form (\ref{gp}) of the
action of $\pa_x$ on a generic element $g$ of the whole Kac-Moody group $G$
makes the derivation easier on $G$.

Let us now say a few words about possible applications and extensions of
this work. First, one may think of the following program of quantization of
our results: to formulate quantum axioms corresponding to nonlocal VPA in
the spirit of \cite{BD}\footnote{In the course of writing this paper, we
became aware of several works dealing with vertex operator structures
containing logarithmic terms, see \cite{Gu,Fl} and references therein}; to
quantize the geometric formulas (\ref{finalPBg-g}), (\ref{finalPBn-n}),
(\ref{localversion}) for nonlocal VPA structures; to realize these formulas
in terms of the algebra of local fields of quantum Toda theories. The
combination of small and large limits for the spectral parameter which we
use is reminiscent of the work \cite{BLZ}. Second, it would be interesting
to carry out the present work in the case of higher AGD structures; this
would lead to a family of compatible nonlocal VPA structures on
$B_{-}\times N_{+}$ (we construct such a family in sect. 4, but its
connection with the AGD structures is not clear to us). Next, it would be
interesting to obtain similar results for other soliton equations such as
the nonlinear Schr\"odinger (NLS) equation; in that case ${\frak a}$ should
be replaced by the loop algebra with values in the Cartan subalgebra
${\frak h}$. The geometric interpretation of the NLS variables analogous to
the one used here, was obtained by Feigin and one of us \cite{FF-toappear}
(the case of $\widehat{\frak{sl}}_2$ was also treated in
\cite{ABF}). Finally, the fact that the VPA structure given in
Thm. \ref{thm3.1} is left $G$-invariant, leads us to conjecture the
existence of affine Weyl group symmetries of the mKdV hierarchies, mixing
local and nonlocal terms (see remark \ref{rem4.2}). These symmetries are
probably connected with the Darboux transformations.

The first author would like to thank B. Feigin for his collaboration
in \cite{EFe} where the ideas of extension of the space of local
fields were used, and A. Orlov and V. Rubtsov for their collaboration
in \cite{EOR} on the subject of nonlocal Poisson structures, and for
discussions about this work. The second author thanks B. Feigin for
useful discussions. He is also grateful to P. Schapira for his
hospitality at Universit\'{e} Paris VI, where this work was completed.

The first author would like to dedicate this paper to A. Guichardet,
his senior colleague who is about to retire from Ecole Polytechnique.

\section{Nonlocal vertex Poisson algebras.}

Let $(R,\pa)$ be a differential ring, and let $\cA$ be the associated
ring of formal pseudodifferential operators. The ring $\cal A$ has
generators $\wt\pa, \wt\pa^{-1}$, and $i(r)$ for $r\in R$, and
relations $\wt \pa\wt\pa^{-1}=\wt\pa^{-1}\wt\pa=1$,
$[\wt\pa,i(r)]=i(\pa r)$ for $r\in R$, and $i:R\to \cal A$ is an
algebra morphism.  In what follows, we will denote by $\pa,\pa^{-1}$
and $r$, $\wt\pa,\wt\pa^{-1}$ and $i(r)$, respectively.

\subsection{The modules $\cR_{n}$.} 

Let $n\ge 1$ be an integer. We will use the following notation in the
algebra $\cA^{\otimes n}$: $\pa_{x_{i}},\pa_{x_{i}}^{-1}$ will stand
for $\otimes_{j=1}^{i-1}1\otimes r\otimes_{j=i+1}^{n}1$,
$\otimes_{j=1}^{i-1}1\otimes \pa\otimes_{j=i+1}^{n}1$,
$\otimes_{j=1}^{i-1}1\otimes \pa^{-1}\otimes_{j=i+1}^{n}1$,
respectively.  For $T\in\cA$, we also denote
$\otimes_{j=1}^{i-1}1\otimes T\otimes_{j=i+1}^{n}1$ by $T(x_{i})$.

We define $\cR_{n}$ to be the quotient of $\cA^{\otimes n}$ by the
left ideal $I_{n}$, generated by $\sum_{i=1}^{n}\pa_{i}$ and
$r(z_{i})-r(z_{j})$ for $r\in R$ and $i,j=1,\ldots,n$. We will
sometimes denote $\cR_{n}$ by $\cR_{n}(R)$. Let us consider $\cR_{n}$
as a left $\cA^{\otimes n}$-module and denote by
$\de_{x_{1},\ldots,x_{n}}$ its generator $1+I_{n}$. We then have the
relations
\begin{equation}
(r(x_{i})-r(x_{j}))\de_{x_{1},\ldots,x_{n}}=0, \quad
(\sum_{i=1}^{n}\pa_{x_{i}})\de_{x_{1},\ldots,x_{n}}=0, 
\label{definedelta}
\end{equation}
so that $\de_{x_{1},\ldots,x_{n}}$ plays the role of the distribution
$\prod_{i=2}^{n}\de(x_{1}-x_{i})$. 

The module $\cR_{2}$ admits the following simple description. 
Let $T\mapsto T^{*}$ 
be the anti-automorphism of $\cA$, defined by $r^{*}=r$
and $\pa^{*}=-\pa$. Let us endow $\cA$ with the $\cA^{\otimes 2}$-module
structure, defined by $(a\otimes b)c=acb^{*}$, for $a,b,c\in\cA$. Then
the linear map $op:\cR_{2}\to\cA$, defined by $op((a\otimes
b)\de_{x_{1}x_{2}})=ab^{*}$, is an isomorphism of $\cA^{\otimes
2}$-modules. The inverse map to $op$ is given by $op^{-1}(a)=(a\otimes
1)\de_{x_{1}x_{2}}= (1\otimes a^{*})\de_{x_{1}x_{2}}$. 

\subsection{Definition of the nonlocal VPA structure.}

A nonlocal VPA structure on $(R,\pa)$ is a linear map $P:R\otimes R\to
\cR_{2}$, satisfying the following conditions:
\begin{equation}
P(ab\otimes c)  =a(x_{1})P(b\otimes c)+b(x_{1})P(a\otimes c),
\label{Leibniz}\end{equation}
\begin{equation}
P(\pa a\otimes b)=\pa_{x_{1}}P(a\otimes b), 
\label{D-linearity}
\end{equation}
\begin{equation}
P(b\otimes a)=-\sigma(P(a\otimes b)),
\label{antisymmetry} 
\end{equation}
for $a,b,c\in R$, where $\sigma$ is the involutive automorphism of
$\cR_{2}$ defined by $\si((a\otimes b)\de_{x_{1}x_{2}})=(b\otimes
a)\de_{x_{1}x_{2}}$, for $a,b\in\cA$, and the Jacobi identity that we
formulate below.

Let us define a map $P_{x,yz}:R\otimes \cR_{2}\to \cR_{3}$ by the
following rules (we attach indices $y,z$ to $\cR_{2}$ and $x,y,z$
to $\cR_{3}$):
$$
P_{x,yz}(a\otimes\de_{y,z})=0, \quad \quad
P_{x,yz}(a\otimes\pa_{y}m)=\pa_{y} P_{x,yz}(a\otimes m),
$$
$$
P_{x,yz}(a\otimes b(y)T(z)\de_{yz})=(op\circ P)(a\otimes
b)(x)T(z)\de_{xyz}+b(y)P_{x,yz}(a\otimes T(z)\de_{yz})
$$
for $a,b\in R$, $T\in\cA$, $m\in\cR_{2}$. 

One can check easily that $P_{x,yz}$ is well-defined by these
conditions. This follows from the identity $P_{x,yz}(a\otimes
\pa_{y}b(y)m)-P_{x,yz}(a\otimes b(y)\pa_{y}m)=P_{x,yz}(a\otimes (\pa
b)(y)m)$, which can be checked by puting $m$ in the form $(1\otimes
T)\de_{yz}$, $T\in \cA$.

The Jacobi identity is then expressed as
\begin{equation}
P_{x,yz}(a\otimes P(b\otimes c))=\si_{xy}[P_{x,yz}(b\otimes
P(a\otimes c))]+\si_{xz}[P_{x,yz}(c\otimes P(b\otimes a))],
\label{Jacobi}
\end{equation}
for any $a,b,c\in R$, where $\si_{xy},\si_{xz}$ are the automorphisms
of $\cR_{3}$, defined by
$$
\si_{xy}(T(x)U(y)V(z)\de_{xyz})=U(x)T(y)V(z)\de_{xyz},
$$
and
$$
\si_{xz}(T(x)U(y)V(z)\de_{xyz})=V(x)U(y)T(z)\de_{xyz},
$$
for $T,U,V\in \cA$.

\medskip
\noindent
\begin{remark} \label{rem1.1}
We may consider elements of $\cR_{2}$ as kernels in $x_{1}$ and
$x_{2}$, expressed as linear combinations of
$r(x_{1})\pa_{x_{1}}^{k}\de(x_{1}-x_{2})$, and think of $P(a\otimes
b)$ as $\{a(x_{1}),b(x_{2})\}$. The terms with $k\ge 0$ are called
local and the terms with $k<0$ are called nonlocal.

The expression $P_{x,yz}(a \otimes m)$ should then be thought of as
expressing the Poisson bracket of the form $\{a(x),m(y,z)\}$, where
$m$ is some kernel.

Note that $P_{x,yz}$ also has the properties
$$
P_{x,yz}(a\otimes\pa_{z}m)=\pa_{z}
P_{x,yz}(a\otimes m),
$$
$$
P_{x,yz}(a\otimes T(y)b(z)\de_{yz})
=(op\circ P)(a\otimes
b)(x)T(z)\de_{xyz}+b(z)P_{x,yz}(a\otimes T(y)\de_{yz}). 
$$

For $a,b,c\in R$, $\{a(x),\{b(y),c(z)\}\}$ is expressed as
$P_{x,yz}(a\otimes P(b\otimes c))$. On the other hand,
$\{b(y),\{a(x),c(z)\}\}$ is expressed as $\si_{xy}[P_{x,yz}(b\otimes
P(a\otimes c))]$, and $\{c(z),\{b(y),$ $a(x)\}\}$ as
$\si_{xz}[P_{x,yz}(c\otimes P(b\otimes a))]$. This explains the
connection between the standard Jacobi identity for the Poisson
brackets and formula (\ref{Jacobi}).
\hfill $\Box$
\end{remark}

\subsection{Connection with the Beilinson-Drinfeld formalism.} 
Let $\cA_{+}$ be the subalgebra of $\cA$, generated by $R$ and $\pa$
(the algebra of differential operators). Let us set
$\cR_{2}^{+}(R)=op^{-1}(\cA_{+})$. We will say that $P$ defines
a local VPA, if $P$ takes values in $\cR_{2}^{+}(R)$.

In this case, the notion described here coincides with that of
``coisson algebra'' (on the disc) of Beilinson and Drinfeld
(\cite{BD}).  Indeed, let $X=\on{Spec}\CC[t]$, and $D_{X}$ be the
ring of differential operators on $X$. Let $A$ be the algebra $R[t]$,
considered as a $D_{X}$-module by the rule that ${d\over {dt}}$ acts
on $R$ as $\pa$. We extend $P$ to an operation
$$
\{,\}\in
\on{Hom}_{D_{X\times X}}
(A\boxtimes A,\Delta_{*}A)
$$ ($\Delta:X\to X \times X$ denotes the diagonal embedding) as in
\cite{BD}, (0.1.7), in the following way.  Let us call $t_{1}=t\otimes
1$ and $t_{2}=1\otimes t$ the coordinates of $X\times X$; $A\boxtimes
A$ is identified with $R\otimes R[t_{1},t_{2}]$, ${\pa\over{\pa
t_{1}}}$ and ${\pa\over{\pa t_{2}}}$ acting on $R\otimes R$ as
$\pa\otimes 1$ and $1\otimes\pa$; on the other hand, $\Delta_{*}A$ is
identified with $\cA_{+}[t]$, with $t_{1,2}$ acting as $t$ and
${\pa\over{\pa t_{1,2}}}$ as ${\pa\over{\pa t}}$.

We then set $\{at_{1}^{n},bt_{2}^{m}\}=t^{n+m}P(a\otimes b)$, for
$a,b\in R$.

\begin{remark} As explained in \cite{BD}, some local VPAs can
be obtained as a classical limit of a family of chiral algebras (or
vertex operator algebras \cite{Bor,FLM}), in the same way as one
obtains Poisson algebras as a classical limit of a family of
associative algebras. In particular, the local VPA $\pi_0$ described
below is the classical limit of the vertex operator algebra of a
Heisenberg algebra (see Remark 2 in \cite{FF-Inv}).

It would be interesting to generalize the notion of chiral algebra to
allow for nonlocality.
\hfill \qed
\end{remark}

\subsection{A class of nonlocal VPA's.}
In this section we give a construction of a class of nonlocal
VPA's. The results of this section will only be used in the proof of
Prp. 3.3. However, the construction presented here might be of general
interest.

\begin{proposition} \label{prop1.1}
Let $E_{k}, k\geq 0$, be the subspace of $\Der(R)^{\otimes 2}$,
consisting of all tensors $\sum_{\al}u_{\al}\otimes v_{\al}$, such that
$$
\sum_{\al}(u_{\al}a)(\pa^{i}v_{\al}b)=0, \quad \forall a,b\in R, 
\quad i=0,\ldots,k.
$$
Suppose we are given elements $\sum_{\al}x_{i}^{(\al)} \otimes
y_{i}^{(\al)}$ of $\Der(R)^{\otimes 2}$ for all $i\geq -1$.

Assume that 
$$
\sum_{\al}x_{i}^{(\al)}\otimes y_{i}^{(\al)}=(-1)^{i+1}\sum_{\al}
y_{i}^{(\al)}\otimes x_{i}^{(\al)}, 
$$
and that 
$$
\sum_{\al}[x_{i}^{(\al)},\pa]\otimes y_{i}^{(\al)}\in \sum_{\al}
x_{i-1}^{(\al)}\otimes y_{i-1}^{(\al)}+E_{i}
$$
(here we set for $i\le -2$, $\sum_{\al} x^{(\al)}_{i}\otimes
y^{(\al)}_{i}=0$; so that $\sum_{\al} x^{(\al)}_{-1}\otimes
y^{(\al)}_{-1}\in (\Der(R)^{\pa})^{\otimes 2}$). 
 
Then the formula
\begin{equation}
P(a\otimes b)=\sum_{i\ge -1}
\sum_{\al}(x_{i}^{(\al)}a)(x)(y_{i}^{(\al)}b)(y) \pa_{x}^{i}\de_{xy},
\label{NLVOA}
\end{equation}
for $a,b\in R$ defines a nonlocal VPA structure on $R$. 
\end{proposition}

{\em Proof. \/} The first condition ensures the antisymmetry of $P$, the
second condition is equivalent to the $\pa$-linearity condition
(\ref{D-linearity}). The first condition being satisfied, the Jacobi 
identity for $P$ is automatically satisfied: for example, the term 
$P_{xy,z}(P(a\otimes b)\otimes c)$ is equal to 
$$
\sum_{i,j,\al,\beta}(x_{j}^{(\beta)}x_{i}^{(\al)}a)(x)(y_{i}^{(\al)}b)(y)
(y_{j}^{(\beta)}c)(z)(-1)^{i+j}\pa_{y}^{i}\pa_{z}^{j}\de_{xyz} 
$$
$$
+\sum_{i,j,\al,\beta}(x_{i}^{(\alpha)}a)(x)(x_{j}^{(\beta)}
y_{i}^{(\al)}b)(y)
(y_{j}^{(\beta)}c)(z)(-1)^{j}\pa_{x}^{i}\pa_{z}^{j}\de_{xyz}
$$
whose second term is cancelled by 
$$
\sum_{i,j,\al,\beta}(x_{j}^{(\beta)}x_{i}^{(\al)}b)(y)(y_{i}^{(\al)}c)(z)
(y_{j}^{(\beta)}a)(x)(-1)^{i+j}\pa_{z}^{i}\pa_{x}^{j}\de_{xyz}
$$
which is a cyclic permutation of the first one. 
\hfill $\Box$

\medskip
\noindent
\begin{remark} \label{rem1.2}
Let us denote by $\Der(R)$ the Lie algebra of derivations of $R$ and
assume that we have $\varpi \in (\Der(R)^{\pa})^{\otimes 2}$, and
$P_{+}:R\otimes R\to \cA_{+}$, such that writing
$\varpi=\sum_{i}\varpi_{i}\otimes \varpi'_{i}$ ($\varpi_{i}$, $\varpi'_{i}$
derivations of $R$ commuting with $\pa$),
\begin{equation}    \label{always}
op\circ P(a\otimes b)=P_{+}(a\otimes b)+\varpi_{i}(a)\pa^{-1}\varpi'_{i}(b). 
\end{equation}
In view of the form of $\varpi$, the nonlocal part of axioms
(\ref{antisymmetry}), (\ref{Jacobi}) is
satisfied. The l.h.s. of the Jacobi identity contains terms
of the form
$$a(x)b(y)c(z)\pa_{x_{i}}^{-1}\pa_{x_{j}}^{-1}, 
a(x)b(y)c(z)\pa_{x_{i}}^{-1}\pa_{x_{j}}^{k},
a(x)b(y)c(z)\pa_{x_{i}}^{k}\pa_{x_{j}}^{\ell},
$$ 
$a,b,c\in R$, $k,\ell\ge 0$, $x_{i}\ne x_{j}$ run through $x,y,z$. The
sum of the terms of the first type then cancels automatically (see
Lemma \ref{lemma2.3}).

All nonlocal VPA structures that we study in this work will be of the
type described here.

Note that in particular, if $P_+$ is $0$, then $P$ given by
(\ref{always}) always defines a VPA structure.

We will denote by $\cA_{-1}$ the span of $\cA_{+}$ and the $a\pa^{-1}
b$, $a,b\in R$, and by $\cR_{2}^{-1}(R)$ the space $op^{-1}(\cA_{-1})$. 
\qed
\end{remark}

\subsection{Hamiltonian vector fields.}

We wish to show briefly here how the notions introduced above can be
related to the Gelfand-Dickey-Dorfman theory of formal variational
calculus (see \cite{GD}, and the introduction of \cite{BD}).  Let us
assume that $P$ takes values in $\cA_{+}$.  Define ${\cal V}_{f}\in\End(R)$
by
\begin{equation}    \label{hamvect}
{\cal V}_{f}(a)=- \left( op \circ P(a\otimes f) \right) \cdot 1,
\end{equation}
the result of the action of the differential operator $P(a\otimes f)$ on
$1\in R$. Then ${\cal V}_{f}$ is a derivation of $R$, commuting with
$\pa$. ${\cal V}_{f}$ actually depends only on the class of $f$ on $R/\pa R$, 
and is called the Hamiltonian vector field corresponding to the density
$f$. We will also use the notation 
$$
{\cal V}_{f}(a)=\{\int_{-\infty}^{\infty} f,a\}=-\{a, \int_{-\infty}^{\infty}
f\}, 
$$
where $\int_{-\infty}^{\infty} f$ denotes the class of $f$ in $R/\pa R$.

It can be used to define a Lie algebra structure on $R/\pa R$, by
the rule, 
$$
\{ \int_{-\infty}^{\infty} f,\int_{-\infty}^{\infty} g \}=
\int_{-\infty}^{\infty} {\cal V}_{f}(g), \quad \quad \forall f,g\in R.
$$

A derivation $D$ of $R$ which commutes with $\pa$
defines an operation $D_{2}$ on $\cR_{2}$, in the following way. We set
$D_{2}\de_{xy}=0$, and extend $D_{2}$ to the whole $\cR_{2}$ by the
condition that it commutes with $\pa_{x}$ and $\pa_{y}$, and the
formulae 
$$
D_{2}(a(x)m(x,y))=a(x)D_{2}m(x,y)+
(Da)(x)m(x,y),
$$
$$
D_{2}(b(y)m(x,y))=b(y)D_{2}m(x,y)+ (Db)(y)m(x,y)
$$
for $m\in
\cR_{2}$, $a,b\in\cA$; this definition makes sense because $D$
commutes with $\pa$.

We then say that $D$ is an infinitesimal automorphism of $R$,
if $P(Da\otimes b)+P(a\otimes Db)=D_{2}P(a\otimes b)$ for $a,b\in R$. 

\begin{proposition} \label{prop1.2}
Under the conditions above, ${\cal V}_{f}$ is an
infinitesimal  
automorphism of $R$. Moreover, $f\mapsto {\cal V}_{f}$ is a Lie algebra
homomorphism from $R/\pa R$ to $\Der(R)^{\pa}$. 
\end{proposition}

{\em Proof\/}. The first part is straightforward and the second is
contained in \cite{GD}. 

\hfill $\Box$

\section{Nonlocal extensions of the VPA of local fields $\pi_{0}$.}

\subsection{Notation and definition of $\pi_{0}$.}

Let $\wt{\frak g}$ be an affine Lie algebra, with generators $e_{i}$,
$f_{i}$, $\al^\vee_{i}$, $i=0,\ldots,l$, and $d$, subject to the
relations of \cite{Kac}. Let $a_{i},a_{i}^{\vee}$ be the labels of
$\widetilde{{\frak g}}$, $h$ be its Coxeter number; then
$K=\sum_{i}a_{i}^{\vee}\al^\vee_{i}$ is a generator of the center of
$\wt{\frak g}$. Let $\widehat{\frak g}$ be the subalgebra of
$\wt{\frak g}$ with the same generators except $d$, ${\frak g}$ be the
quotient of $\widehat{\frak g}$ by $\CC K$, and let us denote the same
way elements of $\wt{\frak g}$ and their images in ${\frak g}$. Let
${\frak h}$ be the Cartan subalgebra of ${\frak g}$, generated by the
$\al^\vee_{i}$'s, ${\frak n}_{+}$ and ${\frak n}_{-}$ be the
pronilpotent subalgebras generated by the $e_{i}$'s and the $f_{i}$'s
respectively, and ${\frak b}_{-}={\frak h}\oplus {\frak n}_{-}$.  Let
$\alpha_{i}$, $i=0,\ldots,l$ be the simple roots of $\wt{\frak g}$,
positive with respect to this decomposition. For each $x \in \g$ we
can write $x=x_+ + x_-$, where $x_+ \in \n_+, x_- \in {\frak b}_-$.

There is an invariant inner product $\langle , \rangle$ on
$\wt{\frak g}$; let us denote in the same way its restriction to
$\widehat{\frak g}$. Let $\sigma$ be any section of ${\frak g}$ to
$\widehat{\frak g}$; the restriction of $\langle , \rangle$ to
$\sigma({\frak g})$ defines an inner product on ${\frak g}$,
independent of $\sigma$ and again denoted by $\langle , \rangle$.

Let $\wt{\frak h}$ be the Cartan subalgebra of $\G$ spanned by
$\al^\vee_i, i=0,\ldots,l$ and $d$. The restriction of the inner
product $\langle , \rangle$ to $\wt{\frak h}$ is non-degenerate and
hence defines an isomorphism $\wt{\frak h} \simeq \wt{\frak h}^*$. Let
$\omega^\vee_i \in \wt{\frak h}^*$ be the $i$th fundamental coweight,
i.e. it satisfies $\langle \omega^\vee_i,\al_j \rangle = \delta_{i,j},
(\omega^\vee_i,d)=0$. Denote by $h_i$, $h^\vee_i$ the elements of
$\wt{\frak h}$, which are the images of $\alpha_i, \omega^\vee_i \in
\wt{\frak h}^*$ under the isomorphism $\wt{\frak h} \simeq \wt{\frak
h}^*$ (note that this is not a standard notation).

Let $$p_{-1}=\sum_{i=0}^{l} \frac{(\al_i,\al_i)}{2} f_{i},$$ and let
${\frak a}$ be the centralizer of $p_{-1}$ in ${\frak g}$. It is a
commutative subalgebra of ${\frak g}$, called the principal abelian
subalgebra. We have ${\frak a}={\frak a}_{+}\oplus{\frak a}_{-}$,
where ${\frak a}_{\pm}={\frak a}\cap {\frak n}_{\pm}$. Let $I$ be the
set (with multiplicities) of integers, congruent to the exponents of
$\wt{\frak g}$ modulo $h$. Then ${\frak a}_{\pm}$ is generated by
elements $p_{n}$, $n\in\pm I$.

We
normalize the $p_{n}$, $n\in \pm I$ is such a way that
$$
\langle p_{n},p_{-n} \rangle ={1\over h}, \quad n\in I, 
$$
where $h$ is the Coxeter number of $\wt{{\frak g}}$.
Let $\pi_{0}$ be the free differential ring generated by $u_{i}$,
$i=1,\ldots,l$; we have $\pi_{0}=\CC[u_{i},\pa u_{i},\ldots]$.

\begin{proposition}[\cite{BD}] \label{prop2.1}
The space $\pi_{0}$ has a VPA structure, defined by the formula
\begin{equation}    \label{PBgen}
P(A\otimes B)= \sum_{1\leq i,j \leq l} (\al_{i},\al_{j}) 
\sum_{n\geq
0} \sum_{m\geq 0} \Big(\frac{\pa A}{\pa u_i^{(n)}} \pa^{n+1} \otimes
\frac{\pa B}{\pa u_j^{(m)}} \pa^m\Big) \delta_{xy}.
\end{equation}
\end{proposition}

\medskip
\noindent
\begin{remark} \label{rem2.1}
This structure is uniquely determined by the following formulas:
$$P(A \otimes 1) = 0, \quad \quad \forall A \in \pi_0,$$ and
\begin{equation}
P( u_{i}\otimes  u_{j})=(\al_{i},\al_{j})\pa_{x}\de_{xy}.
\label{PBphi-phi}
\end{equation}
\hfill $\Box$
\end{remark}

The space $\pi_{0}$ also has the structures of ${\frak n}_{+}$-- and
${\frak a}_{-}$--modules, defined in \cite{FF-CIME}. Below we will
extend these structures in three ways.

\subsection{Extension by half IM's.}
Let $B_{-}$ be the ind-algebraic group corresponding to ${\frak b}_{-}$,
and $N_{+}$ be the pro-algebraic Lie group corresponding to ${\frak
n}_{+}$, $A_{+}$ its subgroup corresponding to ${\frak a}_{+}$. Let $G$ be
the group corresponding to ${\frak g}$, containing $B_{-}$ and $N_{+}$ as
subgroups.

Recall from \cite{FF-CIME}, \cite{FF-Inv} the identification of
$\pi_{0}$ and $\CC[N_{+}/A_{+}]$ as rings and ${\frak
n}_{+}$-modules. Moreover, the Lie algebra ${\frak a}_{-}$ acts on
$\CC[N_+/A_+]$ from the right, since we can identify $N_+$ with an open
subspace of $B_{-}\backslash G$. Let $\pa_n$ be the derivation of
$\pi_0$ corresponding to the right action of $p_{-n}$ on $N_+/A_+$. In
particular, $\pa_1 \equiv \pa$.

Consider for $n\in I$, the hamiltonians
$H_{n}\in \pi_{0}$ from \cite{EFr}. They satisfy
$$
\ep_{-\al_i} e_{i} \cdot H_{n}=\pa A_{n}^{(i)},
$$
for certain $A_{n}^{(i)}\in\pi_{-\al_{i}}$ (in the notation of
\cite{EFr}). We have an isomorphism $\pi_0 \simeq \CC[N_+/A_+]$ (see
\cite{FF-CIME}, \cite{FF-Inv}). 
$$
\pa_{n}H_{m}=\pa_{m}H_{n}=\pa H_{n,m}
$$
for certain $H_{n,m}\in\pi_{0}$. Following \cite{EFr}, Sect. 4, define
$\pi_{0}^{+}=\pi_{0}\otimes \CC[F_{n}]_{n\in I}$, 
and extend the action of $\pa_n$ to it by the formula
\begin{equation}
\pa_{n}F_{m}= H_{n,m}.
\label{actionofhigherflowsonhalfIM}
\end{equation}
(so that $F_{n}$ can be viewed as a ``half integrals of motion''
$\int_{-\infty}^x H_{n}$). In particular, $\pa F_m = H_m$.

According to Thm. 5 and Prop. 9 of \cite{EFr}, the ring $\pi_0^+$ is
isomorphic to $\CC[N_+]$, and the action of $\pa_n$ on $\pi_0^+$
defined this way corresponds to the right action of $p_{-n}$ on $N_+$.

Finally, the action of the generators of ${\frak n}_{+}$ are defined by
\begin{equation}
e_{i} \cdot F_{n}=\ep_{-\al_{i}}^{-1} A_{n}^{(i)} = \phi_n(e_i)
\label{actionofgeneratorsofn-onhalfIM}
\end{equation}
(in the notation of \cite{EFr}).

Now we extend the VPA structure on $\pi_{0}$ to a nonlocal VPA
structure on $\pi_{0}^{+}$. For $a\in \pi_{0}$, $n\in I$, we have
\begin{equation}    \label{inthn}
\{\int_{-\infty}^{\infty}H_{n},a\}=n\pa_{n}a,
\end{equation}
(\cite{DS}, prop. 4.5) so that 
\begin{equation}    \label{onemore}
P(H_{n}\otimes a)\in n\pa_{n}a(y)\de_{xy}+\pa_{x}\cR_{2}^{+}(\pi_{0}).
\end{equation}
Let $i_{n}(a)$ be the element of $\cR_{2}^{+}(\pi_{0})$, such that 
$P(H_{n}\otimes a)= n\pa_{n}a(y)\de_{xy}+\pa_{x}(i_{n}(a))$ 
(this defines
$i_{n}(a)$ uniquely, since $\xi\in\cA$, $\pa\cdot\xi=0$ implies $\xi=0$).
We then define
\begin{equation}
P(F_{n}\otimes a)=n\pa_{n}a(y)
\pa_{x}^{-1}\de_{xy}+i_{n}(a). 
\label{PBF-phi}
\end{equation}
Next, we have:
$$
P(H_{n}\otimes H_{m})=n\pa_{n}H_{m}(x)\de_{xy}
+\pa_{x}(i_{n}(H_{m})). 
$$
Integrating this expression w.r.t. the second variable, we obtain 
$$
n\pa_{n}H_{m}(x)+\pa_{x}(\int_{-\infty}^{\infty}i_{n}(H_{m})(x,y)dy).
$$
where
$$
\int_{-\infty}^{\infty}i_{n}(H_{m})(x,y)dy
$$
is defined as $a_{0}(x)$, where $i_{n}(H_{m}) =
\sum_{i\ge 0}a_{i}(x)\pa_{x}^{i}\de_{xy}$.

On the other hand, we obtain from formula (\ref{inthn}) that
$$
\{H_{n}(x),\int_{-\infty}^{\infty}H_{m}\}=-m\pa_{m}H_{n}(x).
$$

Comparing the last two formulas we find:
$\int_{-\infty}^{\infty}i_{n}(H_{m})(x,y)dy=-(n+m)H_{n,m}(x)$ (for
degree reasons it cannot contain constant terms) and so
$i_{n}(H_{m})=-(n+m)H_{n,m}(x)\de_{xy} +\pa_{y}Q_{n,m}$, with
$Q_{n,m}\in\cA_{+}$. Thus, we obtain:
$$
P(H_{n}\otimes H_{m})=n(\pa_{n}H_{m})(x)\de_{xy}-(n+m)
H_{n,m}(y)\pa_{x}\de_{xy}+(\pa Q_{n,m})(x) \pa_y \delta_{xy}.
$$
Using formula $H_i = \pa F_i$ we finally obtain:
\begin{equation}
P(F_{n}\otimes F_{m})
=(mH_{n,m}(x)+nH_{n,m}(y))\pa_{x}^{-1}\de_{xy}+Q_{n,m}(x)\delta_{xy}.
\label{PBF-F}
\end{equation}

\begin{proposition} \label{prop2.2}
The formulae (\ref{PBphi-phi}), 
(\ref{PBF-phi}), (\ref{PBF-F})
define a nonlocal VPA
structure on $\pi_{0}^{+}$. 
\end{proposition}

{\em Proof\/}. To establish the antisymmetry condition, we should check
that 
$$
P(F_{n}\otimes F_{m})=-\sigma(P(F_{m}\otimes F_{n}));
$$ 
this identity is $\pa_{x}^{-1}\pa_{y}^{-1}$ applied to
the same identity, with $H_{n}$ and $H_m$ in place of
$F_{n}$ and $F_m$, which is true. The same argument works with
the bracket $P(F_{n}\otimes a)$, $a\in\pi_{0}$, with only
one application of $\pa^{-1}$. 

Let us pass to the Jacobi identity. We should check it for the tensors 
$F_{n}\otimes F_{m}\otimes F_{k}$,
$F_{n}\otimes F_{m}\otimes a$, $F_{n}\otimes a\otimes b$,
$a,b\in\pi_{0}$. These identities are 
$\pa_{x}^{-1}\pa_{y}^{-1}\pa_{z}^{-1}$,
resp. $\pa_{x}^{-1}\pa_{y}^{-1}$, $\pa_{x}^{-1}$ applied to the same
identities with
$H$'s replacing the $F$'s, which are true. 
\hfill $\Box$

\medskip

\begin{proposition} \label{prop2.3}
The action of $\pa_{n}$ on $\pi_0^+$ is an infinitesimal automorphism
of the nonlocal VPA structure on $\pi_{0}^{+}$.

Moreover, we have for all $a\in \pi_{0}^{+}$, 
\begin{equation} \label{hamactonpi0+}
P(H_{n}\otimes a)\in n\pa_{n} a(x)\de_{xy} +
\pa_{x}\rho_{n}(a), 
\end{equation}
with $\rho_{n}(a)\in \cR_{2}^{-1}(\pi_{0}^{+})$, such that 
$$
\rho_{n}(a)\in \sum_{m\in
I}H_{n,m}(x)\rho_{n,m}(a)(y)\pa_{x}^{-1}\de_{xy}
+\cR_{2}^{+}(\pi_{0}^{+}), 
$$
where $\rho_{n,m}$'s are linear endomorphisms of $\pi_{0}^{+}$. 
\end{proposition}

{\em Proof \/}. Since $n\pa_{n}$ coincides with ${\cal V}_{H_n}$, it
satisfies the infinitesimal automorphism identity on $\pi_{0}^{\otimes
2}$. To see that this identity is also satisfied for the tensors
$F_{n}\otimes a$ ($a\in\pi_{0}$) and $F_{n}\otimes F_{m}$, we remark
that they are $\pa_{x}^{-1}$, resp. $\pa_{x}^{-1}\pa_{y}^{-1}$ applied
to the similar identities for $H_{n}\otimes a$, resp. $H_{n}\otimes
H_{m}$.

According to formula (\ref{onemore}), equation \ref{hamactonpi0+}) is
satisfied for $a\in \pi_{0}$, with $\rho_{n}(a)$ actually lying in
$\cR_{2}^{+}(\pi_{0})$. Therefore (\ref{hamactonpi0+}) holds for
$a=F_{m}$, as can be seen by applying $\pa_{x}$ to (\ref{PBF-F}). By
Leibnitz rule, if (\ref{hamactonpi0+}) is true for $a, b \in
\pi_0^+$, then it is true for $ab$. Moreover, we find that
$$
\rho_{n}(ab)=a(y)\rho_{n}(b)+b(y)\rho_{n}(a).
$$
This proves (\ref{hamactonpi0+}) and the properties of $\rho_{n}$ in
general.  \hfill $\Box$

\medskip
\noindent
\begin{remark}
Formula (\ref{hamactonpi0+}) can be viewed as a nonlocal substitute of
the hamiltonian property of $H_{n}$. \hfill $\Box$
\end{remark}

\subsection{Extension by $\CC[B_{-}]$.}

Let $\CC[B_{-}]$ be the ring of algebraic functions on $B_{-}$. Let
$\wt \pi_{0}=\pi_{0}\otimes \CC[B_{-}]$. We extend the actions of
${\frak n}_{+}$ and ${\frak a}_{-}$ on $\pi_{0}$, to actions of
${\frak g}$ and ${\frak a}_{-}$ on $\wt\pi_{0}$, in the following way.
 
Recall the identification as rings and ${\frak n}_{+}$-modules, of
$\pi_{0}$ with $\CC[N_{+}/A_{+}]$. Moreover, the ${\frak
a}_{-}$-action on $\pi_{0}$ is identified with its action by vector
fields on $B_{-}\backslash G/A_{+}$ from the right, where
$N_{+}/A_{+}$ is an open subset (recall that $G$ is the group
corresponding to $\frak g$). The ring $\wt\pi_{0}$ is then identified
with the ring of functions on $B_{-}\times N_{+}/A_{+}$.

We define on it the actions of ${\frak g}$ and ${\frak a}_{-}$ as
follows: the value of the vector field generated by $x\in {\frak g}$
at $(b_{-},n_{+}A_{+})$ is
\begin{equation}
(r(\Ad (b_{-}^{-1})(x))_{-}, \ell(\Ad (b_{-}^{-1})(x))_{+})
\end{equation}
and the value of the vector field $\pa_{n}$ (denoted by $\pa$ for
$n=1$) generated by $p_{-n}$, $n\in I$ is
\begin{equation}
(r(\Ad (n_{+})(p_{-n})_{-}), \ell(\Ad (n_{+})(p_{-n})_{+})); 
\end{equation}
here $r(y)$ and $\ell(y)$ denote the right and left vector fields
generated by a Lie algebra element $y$. The proof of these formulas is
analogous to the proof of Lemma 1 in \cite{EFr}.

In particular, for $n=1$ we have according to Lemma 2 of \cite{EFr}
$$
\pa n_{+}= (n_+ p_{-1} n_+^{-1})_+ n_+ = - (n_+ p_{-1} n_+^{-1})_- n_+
+ (n_+ p_{-1} n_+^{-1}) n_+ = 
$$
\begin{equation}
= -(p_{-1}+\sum_{i}u_{i}h_{i}^{\vee})n_{+}+n_{+}p_{-1}.
\label{evoln+}
\end{equation}
We also have:
\begin{equation}
\pa b_{-}= b_- (n_+ p_{-1} n_+^{-1})_- = b_{-}(p_{-1}+\sum_{i}
u_{i}h_{i}^{\vee}).
\label{evolb-}
\end{equation}
The last two formulas should be considered in an arbitrary
representation of $G$ of the form $V((\la))$, where $V$ is
finite-dimensional (see \cite{EFr}).

Set now
\begin{equation}
P(b_{-}\otimes
 u_{i})=-\pa_{y}\big[\Ad(b_{-}(y))(h_{i})b_{-}(x)\pa_{x}^{-1}\de_{xy}\big].
\label{PBb-phi} 
\end{equation}
By this formula we mean the following. In any representation of $B_-$
of the form $V((\lambda))$, an element $b_-$ of $B_-$ can be viewed as
a matrix $(b_{-,kl})$ whose entries $b_{-,kl}$ are Taylor series in
$\la^{-1}$ with coefficients in the ring $\CC[B_-]$. Such functions in
fact generate $\CC[B_-]$. The left hand side of formula
(\ref{PBb-phi}) is the matrix whose entries are $P(b_{-,kl} \otimes
u_i)$. The right hand side of the formula is also a matrix of the same
size, whose entries are elements of $\cR_2(\CC[B_-])$. Via Leibnitz
rule, formula (\ref{PBb-phi}) defines $P(f \otimes u_i) \in
\cR_2(\CC[B_-])$ for any $f \in \CC[B_-]$. We interpret similarly
formulas below for $P(b_- \otimes b_-), P(g \otimes g)$, etc.

\medskip\noindent
\begin{remark} \label{addrem1}
In view of (\ref{evolb-}), we consider $b_{-}(x)$ as the ordered
exponential
$$P\exp\int_{-\infty}^{x}(p_{-1}+\sum_{i=1}^{l}u_{i}(z)h_{i}^{\vee})
dz.$$ Note that in any representation $b_{-}(x)$ is represented by the
matrix
\begin{equation}    \label{scr}
\left( \on{Id} - \sum_{i=0}^l \frac{(\al_i,\al_i)}{2} f_i S_i + \ldots \right)
\exp \left( \sum_{i=1}^l h_i^\vee \varphi_i(x) \right),
\end{equation}
where $\varphi_i(x) = \int_{-\infty}^x u_i(y) dy$, and $S_i$ is the
$i$th ``half screening'' $$S_i = \int_{-\infty}^x e^{-\varphi_i(y)}
dy.$$ The reason for this terminology is that if in the last formula
we integrate over a closed contour, we obtain the classical limit of
the screening operator of conformal field theory; the sum of the
screening operators coincides with the hamiltonian of the affine Toda
field theory. The other terms in the first factor of formula
(\ref{scr}) can be expressed as consecutive Poisson brackets of the
$S_i$'s.

Poisson bracket (\ref{PBb-phi}) can be informally obtained as
follows. Since
$$
P((p_{-1}+\sum_{i=1}^{l}u_{i}h_{j}^{\vee})\otimes \varphi_{j})
=
-h_{j}\de_{xy},
$$
we can write 
$$
\begin{array}{rcl}
\{b_{-}(x),\varphi_{j}(y)\}=
\int_{-\infty}^{x}dz
P\exp\int_{-\infty}^{z}(p_{-1}+\sum_{i=1}^{l}u_{i}h_{i}^{\vee})
(-h_{j}\de_{zx})
\\
P\exp\int_{z}^{x}(p_{-1}+\sum_{i=1}^{l}u_{i}h_{i}^{\vee})
\\
=
-\big(
P\exp\int_{-\infty}^{y}(p_{-1}+\sum_{i=1}^{l}u_{i}h_{i}^{\vee}) \big)
h_{j}
\big( 
P\exp\int_{y}^{x}(p_{-1}+\sum_{i=1}^{l}u_{i}h_{i}^{\vee})
\big)1_{y<x}
\\
=-b_{-}(y)h_{j}b_{-}(y)^{-1}b_{-}(x)1_{y<x},  
\end{array}
$$
where $1_{y<x}$ is the function of $x,y$ equal to $1$ when $y<x$ and
to $0$ else. For $y$ fixed, this is a shifted Heaviside function in
$x$; applying $\pa_{x}$ to it gives $\delta(x-y)$, so that we identify
$1_{y<x}$ with $\pa_{x}^{-1}\de_{xy}$. Formula (\ref{PBb-phi}) is
obtained by applying $\pa_{y}$ to this identity. \hfill $\Box$
\end{remark}

\begin{lemma} \label{lemma2.1}
Definition (\ref{PBb-phi}) is compatible with the $\pa$-linearity
condition (\ref{D-linearity}) on $P$.
\end{lemma}

{\em Proof \/}. Both
left and right hand sides of (\ref{PBb-phi}) satisfy the same identity
$$
\pa_{x}(\on{lhs})=(\on{lhs})(p_{-1}+\sum_{j} u_{j}(x)h_{j}^{\vee})
- b_{-}(x)h_{i}
\pa_{y}\de_{xy},
$$
and
$$
\pa_{x}(\on{rhs})=(\on{rhs})(p_{-1}
+\sum_{j} u_{j}(x)h_{j}^{\vee}) - b_{-}(x)h_{i}
\pa_{y}\de_{xy}.
$$
\hfill $\Box$ 
\medskip

Let us define for $b_{-}\in B_{-}$, and $x\in {\frak g}$,
$r(x)(b_{-})=(\Ad b_{-}(x))_{-}b_{-}$ in any representation of $G$
(recall that $x_{-}$ stands for the projection of $x\in{\frak g}$ on
the second factor of the decomposition ${\frak g}={\frak n}_{+}\oplus
{\frak b}_{-}$). This formula defines the right action of ${\frak g}$
on $B_{-}$ viewed as an open subset of $N_{+}\backslash G$.

Let $\Delta_+$ be the set of positive roots of $\G$, and let
$e^{\al}$, $e_{\al}$, $\alpha \in \Delta_+$, be dual bases of ${\frak
n}_{+}$, ${\frak n}_{-}$ for the inner product $\langle, \rangle$.

Let
\begin{equation}
R^{+}={1\over 2}(\sum_{\al}e^{\al}\otimes
e_{\al}+\sum_{i}h_{i}\otimes h_{i}^{\vee} -\sum_{\al}e_{\al}\otimes
e^{\al}),
\end{equation}
\begin{equation}
R^{-}={1\over 2}(\sum_{\al}e^{\al}\otimes
e_{\al}-\sum_{i}h_{i}\otimes h_{i}^{\vee} -\sum_{\al}e_{\al}\otimes
e^{\al}),
\end{equation}

Let us define, in the tensor product of any pair of representations of
$G$,
\begin{equation}
\begin{array}{rcl}
\everymath{\displaystyle}
P(b_{-}\otimes b_{-}) & = & \{\big((r\otimes r)(R^{-})(b_{-}(x)\otimes
b_{-}(x))\big)(1\otimes b_{-}(x)^{-1}b_{-}(y)) \\
 & & -\big((r\otimes r)(R^{+})(b_{-}(y)\otimes
b_{-}(y))\big)(b_{-}(y)^{-1}b_{-}(x)\otimes 1)
\}
\pa_{x}^{-1}\de_{xy},
\end{array} 
\label{PBb-b}
\end{equation}
where $pr_{-}$ is the projection $\g \rightarrow {\frak b}_-$ along
$\n_+$. This formula can be derived using the argument of Remark
\ref{addrem1}. A similar formula has been obtained by Faddeev and
Takhtajan \cite{Faddeev} (in the rational case) and used by Bazhanov,
Lukyanov and Zamolodchikov \cite{BLZ} (see also \cite{Bab}).

Formula (\ref{PBb-b}) can also be written as
\begin{equation}
\begin{array}{rcl}
\everymath{\displaystyle}
P(b_{-}\otimes b_{-}) =  
(pr_{-}\otimes pr_{-})[\Ad^{\otimes 2}(b_{-}(x))(R^{-})
-\Ad^{\otimes 2}(b_{-}(y))(R^{+})] \\ (b_{-}(x)\otimes b_{-}(y))
\pa_{x}^{-1}\de_{xy}.
\end{array} 
\label{convenientPBb-b}
\end{equation}
It satisfies the antisymmetry condition (\ref{antisymmetry}), because
$$
R^{-}=-R^{+(21)}.
$$

\begin{lemma} \label{lemma2.2}
Formula (\ref{PBb-b}) is compatible with the
$\pa$-linearity condition (\ref{D-linearity}) on $P$. \end{lemma}

{\em Proof\/}. The l.h.s. of (\ref{PBb-b}) satisfies
$$
\begin{array}{rcl}
\everymath{\displaystyle}
\pa_{y}(\on{lhs})=(\on{lhs})(1\otimes (p_{-1}
+\sum_{i} u_{i}(y)h_{i}^{\vee}))
\\
+
\sum_{i}\{b_{-}(x)h_{i}\de_{xy}-\Ad(b_{-}(y))([p_{-1},h_{i}])b_{-}(x)\}
\otimes b_{-}(y)h_{i}^{\vee}\pa_{x}^{-1}\de_{xy},
\end{array} 
$$
and the r.h.s. satisfies
$$
\begin{array}{rcl}
\everymath{\displaystyle}
\pa_{y} (\on{rhs})  = 
\Big\{[(r\otimes r)(R^{-})(b_{-}(x)\otimes b_{-}(x))][1\otimes
b_{-}(x)^{-1}b_{-}(y)(p_{-1} \\  +\sum_{i}h_{i}^{\vee} u_{i}(y))]\\
-\Big[(r\otimes r)(R^{+})
[(b_{-}(y)\otimes b_{-}(y))(p_{-1}\otimes 1+1\otimes
p_{-1} 
-\sum_{i} u_{i}(y)(h_{i}^{\vee}\otimes 1 \\ +1\otimes
h_{i}^{\vee})]\Big]  (b_{-}(y)^{-1}b_{-}(x)\otimes 1) \\
+[(r\otimes r)(R^{+}) (b_{-}(y)\otimes
b_{-}(y))][(p_{-1}+\sum_{i} u_{i}(y)h_{i}^{\vee})b_{-}(y)^{-1}
b_{-}(x)\otimes 1]\Big\}
 \\ \pa_{x}^{-1}\de_{xy} \\
-[(r\otimes r)(R^{-}-R^{+})\big(b_{-}(x)\otimes b_{-}(x)\big)]\de_{xy}
\end{array} 
$$
The $\de_{xy}$--terms coincide because
$R^{+}-R^{-}=\sum_{i}h_{i}\otimes h_{i}^{\vee}$. The
$\pa_{x}^{-1}\de_{xy}$--terms containing $ u_{i}(y)$ are coincide,
since $R^{+}$ is invariant by the conjugation by ${\frak h}$.

The identification of the remaining terms follows from the formula
\begin{equation}    \label{rho}
[p_{-1}\otimes 1 + 1\otimes
p_{-1},R^{+}]=\sum_{i}[p_{-1},h_{i}]\otimes h_{i}^{\vee}.
\end{equation}
Let $\rho:{\frak g}\to {\frak g}$ be the linear map defined by
$\rho(x)=\langle R^{+},1\otimes x\rangle$. Formula (\ref{rho}) means
that $[\ad(p_{-1}),\rho]$ is the linear endomorphism of ${\frak g}$,
equal to $0$ on ${\frak n}_{+}$ and ${\frak n}_{-}$ and to
$\ad(p_{-1})$ on ${\frak h}$, which can be easily verified.
\hfill $\Box$
\medskip

\noindent
\begin{remark} \label{rem2.2}
It follows from the proof that the brackets (\ref{PBb-b}) can also be
expressed replacing $R^{+}$ and $R^{-}$ by
$R^{+}+\kappa\sum_{i}\epsilon^{i}\otimes \epsilon_{i}$ and
$R^{-}+\kappa\sum_{i}\epsilon^{i}\otimes \epsilon_{i}$,
$(\epsilon^{i})$, $(\epsilon_{i})$ dual bases of ${\frak g}$. This
follows also from the fact that $[\sum_{i}\epsilon^{i}\otimes
\epsilon_{i}, b_{-}\otimes b_{-}]=0$ for any $b_{-}\in B_{-}$.  We
prefer the form given here because then (\ref{antisymmetry}) is
manifestly satisfied. \hfill $\Box$
\end{remark}

\medskip
Let us extend $P$ to $\CC[B_{-}]^{\otimes 2}$ by linearity, to
$\CC[B_{-}]\otimes \pi_{0}$ by formulas (\ref{PBb-phi}),
(\ref{D-linearity}) and Leibnitz rule, to $\pi_{0}\otimes \CC[B_{-}]$
by antisymmetry, and finally to $(\CC[B_{-}]\otimes \pi_{0})^{\otimes
2}$ by the Leibnitz rule.

\begin{proposition} \label{prop2.4}
The operation $P$ on $\CC[B_{-}]\otimes \pi_{0}$ defined by formulas
(\ref{PBphi-phi}), (\ref{PBb-phi}), (\ref{PBb-b}) and the rules above,
satisfies the nonlocal VPA axioms.
\end{proposition}

{\em Proof\/}. By construction, $P$ satisfies the $\pa$-linearity
axiom. The Leibnitz rule is satisfied for the brackets involving
$b_{-}$'s since the right hand sides of (\ref{PBb-phi}) and
(\ref{PBb-b}) respect the tensor structure. For the other brackets,
the Leibnitz rule is satisfied by construction. As was already
mentioned, the antisymmetry of $P$ follows from $R^{+}=R^{-(21)}$. Let
us check now the Jacobi identity. For the $\CC[B_-]$ part it follows
from the following lemma.

\begin{lemma}    \label{lemma2.3}
Let $(R,\pa)$ be a differential ring. Let $R_{0}$ be a subring of $R$,
such that $R$ is spanned by elements of the form $\pa^i a_i$, where
$a_i \in R_0$.

Let $\sum_{i}\varpi_{i}\otimes\varpi'_{i}\in S^{2}(\Der(R_{0}))$, 
and let us
assume that there exists a linear map $P:R\otimes R\to\cR_{2}$,
satisfying the Leibnitz rule and the $\pa$-linearity, such that
$$
P(a\otimes b)=\sum_{i}(\varpi_{i}a)(x)(\varpi'_{i}b)(y)
\pa^{-1}_{x}\delta_{xy},
$$
for $a,b\in R_{0}$. Then $P$ satisfies the Jacobi identity.
\end{lemma}
{\em Proof\/}. It is enough to check it for elements of $R_{0}$. 
We have 
$$
\begin{array}{rcl}
\everymath{\displaystyle} P_{xy,z}(P(a\otimes b)\otimes
c)=\sum_{i,j}(\varpi_{j}\varpi_{i}a)(x)
(\varpi'_{i}b)(y)(\varpi'_{j}c)(z)\varpi_{x}^{-1}
\pa_{x}^{-1}\de_{xyz}
\\ +(\varpi_{i}a)(x)(\varpi_{j}\varpi'_{i}b)(y)
(\varpi'_{j}c)(z)\pa_{x}^{-1}\pa_{y}^{-1}\de_{xyz};
\end{array}
$$
the first term is cancelled by 
$$
\sum_{i,j}
(\varpi_{i}c)(z)(\varpi_{j}\varpi'_{i}a)(x)(\varpi'_{j}b)(y)\pa_{z}^{-1}
\pa_{x}^{-1}\delta_{xyz}
$$
(which is obtained from the second by a cyclic permutation).
\hfill $\Box$ \medskip

The Jacobi identity in the case of the proposition is then proved as
follows: introduce the variables $\varphi_{i}$, $i=1,\ldots,l$; let $\wt
R=\CC[B_{-}]\otimes \CC[\varphi_{i}^{(k)}]$, and identify $R$ with a
subalgebra of $\wt R$ by $\pa\varphi_{i}= u_{i}$. Extend the
bracket $P$ to $\wt R$, by the rules 
$$
P(\varphi_{i}\otimes\varphi_{j})=-(\alpha_{i},\alpha_{j})\pa_{x}^{-1}
\delta_{xy},
$$
and 
$$
P(b_{-}\otimes\varphi_{i})
=-\Ad(b_{-}(y))(h_{i})b_{-}(x)\pa_{x}^{-1}\de_{xy}. 
$$
Let now $R_{0}=\CC[B_{-}]\otimes \CC[\varphi_{i}]$; the restriction of $P$
to $R_{0}$ is of the form given in the lemma, so $P$ satisfies the
Jacobi identity. 
\hfill $\Box$

\subsection{Extension by $\CC[B_{-}]$ and half IM's}

Consider now the algebra $\bar\pi_{0}=\pi_{0}\otimes\CC[F_{n}]
\otimes\CC[B_{-}]$. We extend the actions of $\pa$ and of the flows
$\pa_{n}$ on it in the way compatible with the embeddings
$\pi_{0}^{+}\subset\bar\pi_{0}$ and
$\wt\pi_{0}\subset\bar\pi_{0}$. We also define a nonlocal VPA
structure on $\bar\pi_{0}$, in such a way that the above embeddings are
nonlocal VPA morphisms, and we define the brackets $P(F_{n}\otimes
b_{-})$ by 
\begin{equation} \label{PBF-b}
P(F_{n}\otimes b_{-})=\pa_{x}^{-1} P(H_{n}\otimes b_{-});
\end{equation}
we then extend $P$ to $\bar\pi_{0}^{\otimes 2}$ by antisymmetry and the
Leibnitz rule (it is easy to check that (\ref{PBF-b}) is 
compatible with the Leibnitz rule for the products of matrix elements of
$b_{-}$). 

\begin{proposition} \label{prop2.5}
The operation $P$ defined by formulas
(\ref{PBphi-phi}), (\ref{PBF-phi}), (\ref{PBF-F}), (\ref{PBb-phi}),
(\ref{PBb-b}), (\ref{PBF-b}), antisymmetry and the
Leibnitz rule, defines a nonlocal VPA structure on
$\bar\pi_{0}$.
\end{proposition} 

{\em Proof\/}. The same as in Prop. \ref{prop2.1}. 
\hfill $\Box$

\medskip

The rest of this section is devoted to proving that the image of $P$
actually lies in $\cR_{2}^{-1}(\bar\pi_{0})$ (the definition of this
space is in Rem. \ref{rem1.2}). It is enough to prove it for
$P(F_{n}\otimes b_{-})$. 

\begin{lemma} \label{addlemma1}

For certain $A_{ni}\in \cR^{+}_{2}(\pi_{0})$, we have 
\begin{equation} \label{F-u}
P(F_{n}\otimes (p_{-1}+\sum_{i=1}^{l}u_{i}h_{i}^{\vee}))
=\sum_{i=1}^{l}
n(\pa_{n}u_{i})(y)h_{i}^{\vee}\pa_{x}^{-1}\de_{xy}
+A_{ni}h_{i}^{\vee}.
\end{equation}
\end{lemma}

{\em Proof\/}. We have 
$
\{\int_{-\infty}^{\infty} H_{n}, u_{i}\}
=n\pa_{n}(u_{i}),
$
so that 
$
P(H_{n}\otimes u_{i})=n\pa_{n}
u_{i}(y)\de_{xy}+\pa_{x}A_{ni}, 
$
with $A_{ni}\in \cR_{2}^{+}(\pi_{0})$. Now 
$$
P(F_{n}\otimes u_{i})=\pa_{x}^{-1}P(H_{n}\otimes u_{i})=
n(\pa_{n}u_{i}(y))\pa_{x}^{-1}\de_{xy}+A_{ni},
$$
and the statement follows.
\hfill $\Box$

\begin{lemma} \label{addlemma1bis}
$\pa_{n}(\sum_{i=1}^{l}u_{i}(x)h_{i}^{\vee})$ is equal to the Cartan
component of  
$$
-[p_{-1},\Ad(n_{+}(x))(p_{-n})].
$$ 
\end{lemma}

{\em Proof\/}. We will give two proofs of this fact. 

Observe that 
$[p_{-1}, (\log n_{+}(x))_{1}]=-\sum_{i} u_{i}(x)h_{i}^{\vee}$; also
$$
\pa_{n}((\log n_{+})_{1} )=(\Ad (n_{+})(p_{-n}))_{1}
$$
for any $n_{+}\in N_{+}$ (here $\log$
stands for the inverse of the exponential map $\exp: \n_+ \rightarrow
N_+$, which is an isomorphism, and 
the index $1$ means the component of principal degree one of an
element of $\wt{\frak g}$). Hence
$$
\begin{array}{rcl}
\pa_{n}(\sum_{i} u_{i}(x)h_{i}^{\vee})=-\pa_{n}([p_{-1},
(\log n_{+}(x))_{1}]) =-[p_{-1}, \pa_{n}(\log n_{+}(x))_{1}] 
\\ =-[p_{-1},(\Ad
(n_{+}(x))(p_{-n}))_{1}]. 
\end{array}
$$

Alternatively, this statement is a formulation of the Cartan part of
the equation
$$
[\pa+p_{-1}+\sum_{i}u_{i}h_{i}^{\vee},
\pa_{n}-(\Ad(n_{+})(p_{-n}))_{+}]=0, 
$$
which follows from the zero-curvature equation (see e.g. formula (10)
of \cite{EFr}) and the fact that $\Ad(n_{+})(p_{-n})$ commutes with
$\pa+p_{-1}+\sum_{i}u_{i}h_{i}^{\vee}$ (see formula (13) of
\cite{EFr}).
\hfill $\Box$

Let us set $A_{ni}=\sum_{k\ge 0}A_{ni}^{(k)}(y)\pa_{x}^{k}\de_{xy}$,
with $A_{ni}^{(k)}\in \pi_{0}$.

\begin{lemma} \label{addlemma2}
Let $\cZ\in \cR_{2}(\bar\pi_{0})$ be such that 
\begin{equation} \label{diffeqZ}
\pa_{y}\cZ=\cZ (p_{-1}+\sum_{i=1}^{l}u_{i}(y)h_{i}^{\vee})
+b_{-}(y)\big(\sum_{i=1}^{l}
n(\pa_{n}u_{i})(y)h_{i}^{\vee}\pa_{x}^{-1}\de_{xy}
+A_{ni}h_{i}^{\vee}\big).
\end{equation}
Then 
\begin{equation} \label{possibleZ}
\begin{array}{rcl}
\cZ=n[b_{-}(y)(\Ad (n_{+}(y))( p_{-n}))_{-}- \Ad (b_{-}(x))(\Ad (n_{+}(x)) 
(p_{-n}))_{-}b_{-}(y)]\pa_{x}^{-1}\de_{xy} \\
+
\Ad(b_{-}(x))\big(\sum_{i=1}^{l}A_{ni}^{(0)}(x)h_{i}^{\vee}\big)
b_{-}(y)\pa_{y}^{-1}\de_{xy} \\
-\sum_{k\ge 1}\pa_{x}^{k} [\Ad(b_{-}(x))\big(\sum_{i=1}^{l}
A_{ni}^{(k+1)}(x)h_{i}^{\vee}\big) b_{-}(y) \de_{xy}].
\end{array}
\end{equation}
\end{lemma}

{\em Proof\/}. Let 
$$
\cZ_{0}=n[b_{-}(y)(\Ad (n_{+}(y)) (p_{-n}))_{-}
- \Ad (b_{-}(x))(\Ad
(n_{+}(x))(p_{-n}))_{-}b_{-}(y)]\pa_{x}^{-1}\de_{xy}.
$$ 
Then 
$$
\begin{array}{rcl}
\pa_{y}\cZ_{0}-\cZ_{0}(p_{-1}+\sum_{i=1}^{l}u_{i}(y)h_{i}^{\vee})
=n\big\{ b_{-}(y)[p_{-1}+\sum_{i=1}^{l}u_{i}(y)h_{i}^{\vee}
,(\Ad(n_{+}(y))(p_{-n}))_{-}] \\ -b_{-}(y) \big(
[p_{-1}+\sum_{i=1}^{l}u_{i}(y)h_{i}^{\vee}
,\Ad(n_{+}(y))(p_{-n})]\big)_{-} \big\}\pa_{x}^{-1}\de_{xy} \\ =-n
b_{-}(y)
[p_{-1},\big(\Ad(n_{+}(y))(p_{-n})\big)_{1}]\pa_{x}^{-1}\de_{xy} \\ =
n b_{-}(y) \pa_{n}(\sum_{i=1}^{l}u_{i}(y)h_{i}^{\vee})
\pa_{x}^{-1}\de_{xy},
\end{array}
$$
where the last equality follows from Lemma
\ref{addlemma1bis}. We then have 
$$
\pa_{y}(\cZ-\cZ_{0})-(\cZ-\cZ_{0})
(p_{-1}+\sum_{i=1}^{l}u_{i}(y))h_{i}^{\vee})
= b_{-}(y)\sum_{i=1}^{l}A_{ni}h_{i}^{\vee},$$
so that 
$$
\pa_{y}[(\cZ-\cZ_{0})b_{-}(y)^{-1}]=\Ad(b_{-}(y))[\sum_{i=1}^{l}
A_{ni}h_{i}^{\vee}],
$$
and the result. 
\hfill $\Box$

\begin{proposition} \label{addprop1}
\begin{align*}
P(F_{n}\otimes b_{-}) &= n[b_{-}(y)(\Ad (n_{+}(y))( p_{-n}))_{-} \\ &-
\Ad (b_{-}(x))(\Ad (n_{+}(x)) (p_{-n}))_{-}b_{-}(y)] \pa_{x}^{-1}\de_{xy}
\\ &+ \Ad(b_{-}(x))\big(\sum_{i=1}^{l}A_{ni}^{(0)}(x)h_{i}^{\vee}\big)
b_{-}(y)\pa_{y}^{-1}\de_{xy} \\ &-\sum_{k\ge
0}\pa_{x}^{k}[\Ad(b_{-}(x))\big(\sum_{i=1}^{l}
A_{ni}^{(k+1)}(x)h_{i}^{\vee}\big) b_{-}(y) \de_{xy}]
\end{align*}
and so belongs to $\cR_{2}^{-1}(\wt{\pi}_{0})$. 
\end{proposition}

{\em Proof\/}.  Indeed, $P(F_{n}\otimes b_{-})$ defined by formula
(\ref{PBF-b}) satisfies (\ref{diffeqZ}), by the Leibnitz rule and
Lemma \ref{addlemma1}, and we apply to it Lemma \ref{addlemma2}.
\hfill $\Box$

\begin{remark} \label{rem2.3}
It follows from Thm. \ref{thm3.1} that the extension of
$\pa_{n}$ to $\bar{\pi}_{0}$, defined again as $r(p_{-n})$, is an
infinitesimal automorphism of the nonlocal VPA structure of
$\bar\pi_{0}$.\hfill $\Box$
\end{remark}

\section{Geometric interpretation of the Poisson structures.}
\subsection{Determination of the nonlocal terms.}
In 2.3, we have defined a nonlocal VPA $\bar\pi_{0}$. It is
isomorphic as a differential algebra to $\CC[B_{-}\times N_{+}]$
with the derivation $\pa$ defined by the right action of $p_{-1}$ on 
$B_{-}\times N_{+}$. We will now describe $P$ in these geometric
terms.  

Let us set $t=\sum_{\al}e^{\al}\otimes e_{\al}+e_{\al}\otimes
e^{\al}+\sum_{i}h_{i}\otimes h_{i}^{\vee}$. Recall that $p_{-n}$
is a basis of ${\frak a}_{-}$ dual to $p_{n}$ with respect to the inner
product $\langle , \rangle$.
 
\begin{lemma}   \label{lemma3.1}
\begin{equation}
P(b_{-}\otimes b_{-})=-\ell^{\otimes 2}(t)(b_{-}(x)\otimes b_{-}(y))
\pa_{x}^{-1}\de_{xy}
\label{reformulationPBb-b}
\end{equation}
\end{lemma}

{\em Proof \/}. 
We have 
$$
\begin{array}{lcl}
\everymath{\displaystyle}
(pr_{-}\otimes pr_{-})(\Ad^{\otimes 2}b_{-}(x)R^{-} - 
\Ad^{\otimes 2}(b_{-}(y))(R^{+})
) =  (pr_{-}\otimes pr_{-})(\Ad^{\otimes 2}(b_{-}(x)) \\ \sum_{\al}e^{\al}
\otimes e_{\al}) 
+(pr_{-}\otimes pr_{-})(\Ad^{\otimes 2}(b_{-}(y))\sum_{\al}e_{\al}
\otimes e^{\al}) -\sum_{i}h_{i}\otimes h_{i}^{\vee}
\\
= \sum_{\al}(\Ad (b_{-}(x))(e^{\al}))_{-}\otimes \Ad (b_{-}(x))(e_{\al})
+\sum_{\al} \Ad (b_{-}(y))(e_{\al}) \otimes (\Ad (b_{-}(y))(e^{\al}))_{-}
\\ -\sum_{i}h_{i}\otimes h_{i}^{\vee}
\\ =-\Ad (b_{-}(x))
(\Ad (b_{-}(x)^{-1})(e^{\al}))_{-}\otimes e_{\al}
-e_{\al}\otimes \Ad (b_{-}(y))(\Ad
(b_{-}(y)^{-1})(e^{\al}))_{-}
\\ -\sum_{i}h_{i}\otimes h_{i}^{\vee}. \end{array}
$$

The first equality is obtained using the following arguments: 
$R^{-}=\sum_{\al}e^{\al}\otimes
e_{\al}-{1\over 2}t$, 
$R^{+}=-\sum_{\al}e_{\al}\otimes e^{\al}+{1\over 2}t$;
$t$ is $B_{-}$-invariant, and $(pr_{-}\otimes pr_{-})t
=\sum_{i}h_{i}\otimes
h_{i}^{\vee}$; the second equality is straightforward, and the 
the last one is because for $b_{-}\in B_{-}$, we have 
\begin{equation}
\sum_{\al}(\Ad (b_{-})(e^{\al}))_{-}\otimes \Ad (b_{-})(e_{\al})
=-\sum_{\al}\Ad (b_{-})(\Ad (b_{-}^{-1})e^{\al})_{-}\otimes e_{\al}.
\label{auxil}
\end{equation}
Formula (\ref{auxil}) can be shown as follows: let $\xi\in{\frak
n}_{+}$, then
$$
\langle \on{lhs\ of\ (\ref{auxil})}, 1\otimes\xi\rangle=(\Ad b_{-}(\Ad
b_{-}^{-1}(\xi))_{+})_{-},
$$ 
and 
$$
\langle \on{rhs\ of\ (\ref{auxil})}, 1\otimes\xi\rangle=-\Ad b_{-}(\Ad
b_{-}^{-1}(\xi))_{-},
$$ 
with $x_{+}=x-x_{-}$. 

Formula (\ref{reformulationPBb-b}) then follows from
(\ref{convenientPBb-b}).  \hfill $\Box$ \medskip

\begin{lemma} [\cite{EFr}, Prop. 6] \label{lemma3.3} 
In any representation of $G$, $n_{+}(x)$ has the form
\begin{equation}
\bar n_{+}(x) \exp \left( \sum_{n\in I} {1\over n}{{p_{n}F_{n}(x)}}
\right),
\label{formofn+}
\end{equation}
where $\bar n_{+}(x)$ is a matrix of polynomials with entries in
$\CC[u_i^{(n)}]$.
\end{lemma}

Recall from \cite{EFr}, Lemma 1, that in any representation of $N_+$
we have the following formula for the right action of $x \in \g$ on
$N_{+}\subset B_{-}\backslash G$: $$r(x)n_{+}=(\Ad
(n_{+})(x))_{+}n_{+}.$$

\begin{lemma} \label{lemma3.4}
\begin{equation}
\begin{array}{rcl}
P(n_{+}\otimes n_{+})= r(a)[n_{+}(x)\otimes
n_{+}(y)]\pa_{x}^{-1}\de_{xy}
+\on{\ local\  terms,}
\end{array}
\label{PBn-n}
\end{equation}
where $a=\sum_{n\in I\cup(-I)}p_{n}\otimes p_{-n}$. 
\end{lemma}

{\em Proof \/}.  Let us write $n_{+}$ in the form (\ref{formofn+}) and
compute the nonlocal part of $P(n_{+}\otimes n_{+})$. It comes from
three different terms: $P(F_{n}\otimes \bar n_{+}), P(\bar n_{+}
\otimes F_n), P(F_{n}\otimes F_m)$.

We have 
$$
P(F_{n}\otimes n_{+})=
n\pa_{n}n_{+}(y)\pa_{x}^{-1}\de_{xy}+\rho_{n}(n_{+}), 
$$
by (\ref{hamactonpi0+}). Set $F=\sum_{n\in I}{1\over
n}p_{n}F_{n}$. We then have $\bar n_{+}=n_{+}e^{-F}$. 
The Leibnitz rule gives 
\begin{align*}
P(F_{n}\otimes \bar n_{+})  =n\pa_{n}n_{+}(y)e^{-F(y)}\pa_{x}^{-1}\de_{xy}
& +\rho_{n}(n_{+})(y)e^{-F(y)}
+(1\otimes \bar n_{+}(y)) \\ 
& [\sum_{m\in I}
{-}{1\over
m}(1\otimes p_{m})P(F_{n}\otimes F_{m})]. 
\end{align*}
Using (\ref{PBF-F}), we obtain 
\begin{align*}
P(F_{n}\otimes \bar n_{+})  & \in 
n\pa_{n}n_{+}(y)e^{-F(y)}\pa_{x}^{-1}\de_{xy}
+\rho_{n}(n_{+})e^{-F(y)}
+(1\otimes \bar n_{+}(y)) \cdot
\\ 
& \left( \sum_{m\in I}{-}{1\over
m}(1\otimes p_{m})(mH_{n,m}(x)+n H_{n,m}(y))\pa_{x}^{-1} \right) \de_{xy}
+\cR_{2}^{+}(\pi_{0}),
\end{align*}
so that 
\begin{align} \label{tmp}
P(F_{n}\otimes \bar n_{+}) &  \in n[(\Ad \bar n_{+}(p_{-n}))_{+}\bar
n_{+}
-\bar n_{+}\sum_{m\in I}{1\over m}p_{m}H_{n,m}](y)\pa_{x}^{-1}\de_{xy}
\\ 
& 
+\rho_{n}(n_{+})e^{-F(y)} 
-\sum_{m\in I}H_{n,m}(x)(\bar
n_{+}(y)p_{m})\pa_{x}^{-1}\de_{xy}+\cR_{2}^{+}(\pi_{0}). 
\end{align}
Let us show that 
$$
\cW=\rho_{n}(n_{+})e^{-F(y)}-\sum_{m\in I}H_{n,m}(x)(\bar
n_{+}(y)p_{m})\pa_{x}^{-1}\de_{xy} 
$$ 
belongs to $\cR_{2}^{+}(\pi_{0})$. Apply $\pa_{x}$ to (\ref{tmp}). The
l.h.s. of the resulting identity belongs to $\cR_{2}^{+}(\pi_{0})$, as
well as 
$
\pa_{x}
(n[(\Ad \bar n_{+}(p_{-n}))_{+}\bar
n_{+}
-\bar n_{+}\sum_{m\in I}{1\over
m}p_{m}H_{n,m}](y)\pa_{x}^{-1})\de_{xy}.
$  
It follows that 
$$
\pa_{x}\cW\in \cR_{2}^{}{+}(\pi_{0}). 
$$
On the other hand, in view of Prop. \ref{prop2.3}, we can write the
nonlocal part of $\cW$ as $\sum_{m\in
I}H_{n,m}(x)w_{m}(y)\pa_{x}^{-1}\de_{xy}$. The result of the action of
$\pa_{x}$ on this nonlocal part should be local; since all $\pa
H_{n,m}$, $m\in I$, are independent, it follows that all $w_{m}$'s
vanish, so that 
$$
\cW\in \cR_{2}^{+}(\pi_{0}).
$$ 

Formula (\ref{tmp}) becomes 
$$
P(F_{n}\otimes \bar n_{+})\in n \left( (\Ad \bar
n_{+}(p_{-n}))_{+}\bar n_{+} -\bar n_{+}\sum_{m\in I}{1\over
m}p_{m}H_{n,m} \right)(y)\pa_{x}^{-1}\de_{xy} +\cR_{2}^{+}(\pi_{0}).
$$
We derive from this a similar statement on $P(\bar n_{+}\otimes
F_{n})$. Combining these with 
(\ref{PBF-F}), and the fact that $P(\bar
n_{+}\otimes \bar n_{+})$ is contained in $\cR_{2}^{+}(\pi_{0})$, 
we obtain (\ref{PBn-n}). 
\hfill $\Box$
\medskip

Introduce the following notation. For any tensor $\gamma = \sum_i
\gamma_i \otimes \gamma'_i$, and two operators $a$ and $b$, we write
$a (\gamma^{(1)}) \otimes b(\gamma^{(2)})$ for $\sum_i a(\gamma_i)
\otimes b(\gamma'_i)$.

\begin{lemma} \label{lemma3.5}
\begin{equation}
\begin{array}{rcl}
P(b_{-}\otimes n_{+})=-b_{-}(x)[\Ad (b_{-}^{-1}(x))(t^{(1)})]_{-}\otimes
[\Ad (b_{-}^{-1}(y))(t^{(2)})]_{+}n_{+}(y)\pa_{x}^{-1}\de_{xy} \\
+\sum_{n\in I}b_{-}(x)(\Ad (n_{+}(x))(p_{-n}))_{-}\otimes
n_{+}(y)p_{n}\pa_{x}^{-1}\delta_{xy}
+\on{\ local\  terms;}
\end{array}
\label{prePBb-n}
\end{equation}
this equation should be understood in the tensor product of
two representations of ${\frak g}$. 
\end{lemma}

The proof is given in Sect. 5.

\subsection{Determination of VPA structures}

\begin{theorem} \label{thm3.1}
The nonlocal VPA structure of $\bar \pi_{0}$
is expressed, via the identification of $\bar\pi_{0}$ with
$\CC[B_{-}\times N_{+}]$,
by the formula 
\begin{equation}
\begin{array}{rcl}
P(g\otimes g)=[-\ell^{\otimes 2}(t)(g(x)\otimes g(y))+r^{\otimes
2}(a)(g(x)\otimes g(y))]\pa_{x}^{-1}\delta_{xy} \\ +
\sum_{n\ge 0}r^{\otimes 2}\Big((\ad p_{-1})^{-n-1}\otimes
1)(t-a)\Big)(g(x)\otimes g(y))\pa_{x}^{n}\de_{xy}.
\label{finalPBg-g}
\end{array}
\end{equation}
\end{theorem}

{\em Proof \/}. Let us denote by $g(x)$ the pair
$(b_{-}(x),n_{+}(x))$; it lies in the variety $B_{-}\times
N_{+}$. This variety is endowed with left and right actions of ${\frak
g}$, that we denote $\ell$ and $r$. They are defined as follows.  The
mapping $B_{-}\times N_{+}\to G$, associating to $(b_{-},n_{+})$ the
product $b_{-}n_{+}$ embeds $B_{-}\times N_{+}$ in $G$ as a Schubert
cell. The left and right actions of ${\frak g}$ on $G$ can be
restricted to $B_{-}\times N_{+}$; so that $\ell(x)$ is the vector
field equal at $(b_{-},n_{+})$ to the sum of $r((\Ad b_{-}(x))_{-})$
on the first component of the product, and $\ell((\Ad b_{-}(x))_{+})$
on the second one, according to
$$
x \cdot b_{-}n_{+}=b_{-}[(\Ad (b_{-})(x))_{-}+(\Ad (b_{-})(x))_{+}]n_{+}.
$$
Likewise, $r(x)$ is the vector field equal at $(b_{-},n_{+})$ to the 
sum
of $r((\Ad n_{+}(x))_{-})$ on the first component of the product, and
$\ell((\Ad n_{+}(x))_{+})$ on the second one, since 
$$
b_{-}n_{+}\cdot x=b_{-}[(\Ad (n_{+})(x))_{-}+(\Ad (n_{+})(x))_{+}]n_{+}.
$$

Formulas (\ref{PBb-b}) and (\ref{PBb-n}) then imply
\begin{equation}
P(g\otimes g)=[-\ell^{\otimes 2}(t)(g(x)\otimes g(y))+r^{\otimes
2}(a)(g(x)\otimes g(y))]\pa_{x}^{-1}\delta_{xy}+\on{ \ local \ terms.}
\label{PBg-g}
\end{equation}

The action of $\pa$ on $g$ coincides with $r(p_{-1})$, due to
formulas (\ref{evolb-}) and (\ref{evoln+}), and we obtain
\begin{equation}
\pa g=r(p_{-1})g.
\label{evolg}
\end{equation}
Therefore $P(g\otimes g)$ satisfies the differential equations
$$
\pa_{x}P(g\otimes g)=(r(p_{-1})\otimes 1)P(g\otimes g), \quad
\pa_{y}P(g\otimes g)=(1\otimes r(p_{-1}))P(g\otimes g). 
$$

We will use these equations to determine $P(g \otimes g)$ completely.
Let us first determine the local part of the r.h.s. of
(\ref{PBg-g}). Denote it by $\cY$.

We have 
\begin{equation}
\pa_{x}\cY- (r(p_{-1})\otimes 1)\cY
=r^{\otimes 2}(t-a)(g(x)\otimes g(x))\de_{xy},
\label{eqlocal1}
\end{equation}
\begin{equation}
\pa_{y}\cY-(1\otimes r(p_{-1}))\cY
=-r^{\otimes 2}(t-a)(g(x)\otimes g(x))\de_{xy}
\label{eqlocal2}
\end{equation}
(we have used the $N_{+}$- and $B_{-}$-invariances of $t$ to replace
$\ell^{\otimes 2}(t)$ by $r^{\otimes 2}(t)$, and that
$(g(x)\otimes g(y))\de_{xy}=(g(x)\otimes g(x))\de_{xy}$).

Recall that ${\frak g}={\frak a}\oplus\Imm (\ad p_{-1})$; this is an
orthogonal decomposition in ${\frak g}$; $\ad p_{-1}$ is
an automorphism of $\Imm(\ad p_{-1})$, and $t-a$ belongs to 
$(\Imm(\ad
p_{-1}))^{\otimes 2}$. 

Then 
\begin{equation}
\cY_{0}= \sum_{n\ge 0}r^{\otimes 2}\Big((\ad p_{-1})^{-n-1}\otimes
1)(t-a)\Big)(g(x)\otimes g(y))\pa_{x}^{n}\de_{xy}
\label{Azero}
\end{equation}
is a solution to (\ref{eqlocal1}); to see it, one should apply the
Leibnitz rule, formula (\ref{evolg}) and note the cancellation of all
terms except the one in $\de_{xy}$. It is also a solution of
(\ref{eqlocal2}) because $((\ad p_{-1})^{n}\otimes
1)(t-a)=(-1)^{n}(1\otimes (\ad p_{-1})^{n})(t-a)$. Indeed, this
identity can be proved by writing $t-a=\sum_{\al}\widehat
e^{\al}\otimes\widehat e_{\al}$, where $\widehat e^{\al}$, $\widehat
e_{\al}$ are dual bases of $\Imm(\ad p_{-1})$, and using the
anti-selfadjointness of $\ad p_{-1}$.

Let $P_{0}$ be the operation defined by (\ref{PBg-g}) and
(\ref{Azero}). It defines a nonlocal VPA structure on
$(\CC[B_{-}\times N_{+}],r(p_{-1}))$ by virtue of Prop. \ref{prop1.1} 
(the elements of $E_{i}$ in the second condition being $\ell^{\otimes
2}(t)-r^{\otimes 2}(t)$
for $i=0$, and $0$ for $i>0$).

The theorem follows from the following lemma which is proved in Sect. 5.

\begin{lemma} \label{lemma3.6}
$P_{0}$ is equal to $P$. 
\end{lemma}
\hfill \qed

\begin{remark} \label{rem3.1}
The variety $B_{-}\times N_{+}$ can be considered as an open subset of
$G$. One can show that formula (\ref{compatPBg-g}) defines a nonlocal
VPA structure on the whole group $G$.
\hfill \qed
\end{remark}

It is easy to derive now a formula for the nonlocal  
VPA structure of $\CC[N_{+}]$.

\begin{corollary} \label{corollary3.1}
The nonlocal VPA structure of $\pi_{0}^{+}$
is expressed, via the identification of $\pi_{0}^{+}$ with $\CC[N_{+}]$,
by the formula 
\begin{equation}
\begin{array}{rcl}
P(n_{+}\otimes n_{+})=r^{\otimes
2}(a)(n_{+}(x)\otimes n_{+}(y))\pa_{x}^{-1}\delta_{xy} \\ +
\sum_{n\ge 0}r^{\otimes 2}\Big((\ad p_{-1})^{-n-1}\otimes
1)(t-a)\Big)(n_{+}(x)\otimes n_{+}(y))\pa_{x}^{n}\de_{xy}.
\label{finalPBn-n}
\end{array}
\end{equation}
\end{corollary}

We now reformulate (\ref{finalPBn-n}) so as to make it clear that the
$A_+$--invariant functions of $N_{+}$ only have local Poisson brackets
expressed in terms of $A_+$--invariant functions of $N_{+}$. Let us
denote
$$
(t-a)_{k}=\big((\ad p_{-1})^{-k-1}\otimes 1\big)(t-a).
$$

\begin{corollary} \label{corollary3.2}
\begin{align} \label{nbarversion}
& P(n_{+}\otimes n_{+})=r^{\otimes
2}(a)(n_{+}(x)\otimes n_{+}(y))\pa_{x}^{-1}\delta_{xy} \\ \notag & +
\sum_{k\ge 0}
\big\{\Ad (\bar n_{+}(x))\big[(\pa_{x}-\sum_{n\in I}{1\over n}H_{n}(x)
\ad p_{n})^{k}((t-a)_{k}^{(1)}\big]
\big\}_{+} n_{+}(x)
\\ \notag
& \otimes
\big(
\Ad(\bar n_{+}(y))(t-a)_{k}^{(2)}
\big)_{+}n_{+}(y) \de_{xy}.
\end{align}
\end{corollary}

{\em Proof \/}. Let $F(x)=\sum_{n\in I}{1\over n}p_{n}F_{n}(x)$, and
$\alpha \in\G$. Then
\begin{equation} \label{identforlocal}
\Ad (e^{F(x)-F(y)})(\alpha\pa_{x}^{n}\de_{xy})=[\pa_{x}-\ad
(F'(x))]^{n}(\alpha\de_{xy}). 
\end{equation}
Indeed, this is obvious for $n=0$. Assume that it is true for $n$, and
apply $\pa_{x}$ to the corresponding identity. We find
$$
\begin{array}{rcl}
\Ad(e^{F(x)-F(y)})(\alpha\pa_{x}^{n+1}\de_{xy})
+[F'(x), \Ad(e^{F(x)-F(y)})(\alpha\pa_{x}^{n}\de_{xy})]
=
\\
\pa_{x}\{[\pa_{x}-\ad
(F'(x))]^{n}(\alpha\de_{xy})\}; 
\end{array}
$$
but the l.h.s. of this identity is expressed as 
$$
\Ad(e^{F(x)-F(y)})(\alpha\pa_{x}^{n+1}\de_{xy})
+[F'(x), \big(\pa_{x}-\ad
(F'(x))\big)^{n}(\alpha\de_{xy})],
$$
and we obtain (\ref{identforlocal}) at step $n+1$.

Now formula (\ref{nbarversion}) follows directly from
(\ref{finalPBn-n}) and (\ref{identforlocal}). In fact, the local part
of formula (\ref{finalPBn-n}) can be rewritten as
$$\sum_{k\geq 0} \left( \Ad(n_+(x) \otimes n_+(y)) (t-a)_k\right)_+
n_+(x) \otimes n_+(y) \pa_x^k \de_{xy}.$$ After substituting
$n_+(x) = \bar{n}_+(x) e^{F(x)}$ in this formula we obtain
$$\sum_{k\geq 0} \left( \Ad(\bar{n}_+(x) \otimes \bar{n}_+(y))
\Ad(e^{F(x)} \otimes e^{F(y)}) (t-a)_k\right)_+ n_+(x) \otimes n_+(y)
\pa_x^k \de_{xy},$$ which is equal to
$$\sum_{k\geq 0} \left( \Ad(\bar{n}_+(x) \otimes \bar{n}_+(y))
(\Ad(e^{F(x)-F(y)} \otimes 1) (t-a)_k\right)_+ n_+(x) \otimes n_+(y)
\pa_x^k \de_{xy}$$ by ${\frak a}_+$--invariance of $(t-a)_k$. Applying
formula (\ref{identforlocal}) we obtain (\ref{nbarversion}).
\hfill
$\Box$

Now let $V((\la))$, $W((\mu))$ be $\G$-modules, and $v\in V$, $w\in W$
be such that ${\frak a_{+}}v=0$, ${\frak a_{+}}w=0$. The matrix
coefficients of $n_{+}(x)v$, $n_{+}(y)w$ are composed of functions on
$N_{+}/A_{+}$, or, equivalently, of elements of $\pi_{0}$. Moreover,
all elements of $\CC[N_{+}/A_{+}]$ can be obtained this way. Formula
(\ref{nbarversion}) implies that
\begin{align} \label{localversion}
& P(n_{+}v\otimes n_{+}w)= \\ \notag
& \sum_{k\ge 0}
\big\{\Ad (\bar n_{+}(x))\big[(\pa_{x}-\sum_{n\in I}{1\over n}H_{n}(x)
\ad p_{n})^{k}((t-a)_{k}^{(1)}\big]
\big\}_{+} n_{+}(x)v \\ \notag & \otimes
\big( \Ad(\bar n_{+}(y))(t-a)_{k}^{(2)} \big)_{+}n_{+}(y)w \de_{xy}; 
\end{align}
since $\bar n_{+}(x)$, $\bar n_{+}(y)$ are matrices whose entries are
elements of $\pi_{0}$, (\ref{localversion}) shows at the same time that
the Poisson brackets of the entries of $n_{+}(x)v$ and $n_{+}(y)w$ are
local and expressed in terms of functions of $N_{+}/A_{+}$. 

In the theory of the mKdV equations, an important role is played by
the embedding of $N_{+}/A_{+}$ in $\G$ as an $N_+$--coadjoint orbit of
an element of ${\frak a}$.  This motivates us to compute the Poisson
brackets of the matrices $\Ad(n_{+})(p_{n})$. The result is a
direct consequence of Cor. \ref{corollary3.1}.

\begin{corollary} \label{corollary3.3} For $n,m\in\pm I$, 
\begin{align*}
& P(\Ad(n_{+})(p_{n})\otimes \Ad(n_{+})(p_{m}))=
\\
& \sum_{k\ge 0}
\big[
\big\{\Ad (\bar n_{+}(x))\big[(\pa_{x}-\sum_{n\in I}{1\over n}H_{n}(x)
\ad p_{n})^{k}((t-a)_{k}^{(1)})\big]
\big\}_{+}, 
\Ad(n_{+}(x))(p_{n})\big]
\\
& \otimes
\big[
\big(
\Ad(\bar n_{+}(y))(t-a)_{k}^{(2)}
\big)_{+}, \Ad(n_{+}(y))(p_{m})
\big]\de_{xy}. 
\end{align*}
\end{corollary}

Here again, the locality of these brackets and the fact they are expressed
in terms of $A_+$--invariant functions on $N_{+}$ is manifest.

\medskip
\noindent
\begin{remark}
Formulas (\ref{finalPBn-n}) and (\ref{localversion}) give a geometric
interpretation of the nonlocal VPA structures of $\pi_{0}^{+}$ and
$\pi_{0}$ respectively.  \hfill $\Box$
\end{remark}
\medskip

We now give a formula for $P(u_{i}\otimes n_{+})$. Due to the presence
of nonlocal quantities $F_{n}$ in $n_{+}$, this bracket contains
nonlocal as well as local terms.

\begin{proposition} \label{prop3.1}
$$
\begin{array}{rcl}
P(u_{i}\otimes n_{+})=n_{+}(y)\pa_{x}\Big(-(n_{+}(x)^{-1}h_{i}
n_{+}(x))_{{\frak
a}}\pa_{y}^{-1}\de_{xy} \\ +\sum_{n\ge 0}(-\ad\
p_{-1})^{-n-1}((n_{+}(x)^{-1}h_{i}n_{+}(x))_{\Imm(\ad p_{-1})})
\pa_{x}^{n}\de_{xy}\Big), 
\end{array}
$$
where the indices ${\frak a}$ and $\Imm(\ad p_{-1})$ stand for the
projections on the components of ${\frak g}={\frak a}\oplus\Imm(\ad
p_{-1})$.
\end{proposition}

{\em Proof. \/}
Consider again the variable
$\varphi_{i}=\pa^{-1}u_{i}$. It is easier to compute
$\cG=P(\varphi_{i}\otimes n_{+})$ first. We have
$$
\pa_{y}\cG=-(p_{-1}+\sum_{i}u_{i}(y)h_{i}^{\vee})\cG+\cG p_{-1}
-h_{i}n_{+}(x)\pa_{xy}. 
$$
Let $\cH=n_{+}(y)^{-1}\cG$, then
$$
\pa_{y}\cH=[\cH,p_{-1}]-n_{+}(x)^{-1}h_{i}n_{+}(x)\pa_{xy}.
$$
A solution to this equation is 
$$
\begin{array}{rcl}
\cH_{0}=-(n_{+}(x)^{-1}h_{i}n_{+}(x))_{{\frak
a}}\pa_{y}^{-1}\de_{xy} \\ +\sum_{n\ge 0}(-\ad\
p_{-1})^{-n-1}((n_{+}(x)^{-1}h_{i}n_{+}(x))_{\Imm(\ad p_{-1})})
\pa_{x}^{n}\de_{xy}.
\end{array}
$$
On the other hand, the nonlocal part of $\cH$ coincides with that of
$\cH_{0}$ because it is equal to the nonlocal part of
$\sum_{n\in I}{1\over n}p_{n}P(\varphi_{i}\otimes F_{n})$; but
$P(\varphi_{i}\otimes F_{n})=-n\pa_{n}\varphi_{i}(x)+\on{local\ terms}$,
since $P(\varphi_{i}\otimes H_{n})\in -n\pa_{n}\varphi_{i}(x)
+\pa_{y}(\cR_{2}(\pi_{0}))$; so this nonlocal part is expressed as
$-\sum_{i\in I}\pa_{n}\varphi_{i}(x)p_{n}\pa_{y}^{-1}\de_{xy}$. To
establish 
the coincidence of the nonlocal parts of $\cH$ and $\cH_{0}$, it remains to
check the identity 
\begin{equation}
\pa_{n}\varphi_{i}(x)=\langle n_{+}^{-1}(x)h_{i}n_{+}(x), p_{-n}\rangle
\label{evolphi}
\end{equation}
which amounts to 
$$
\pa_{n}u_{i}(x)=-\langle h_{i}, [p_{-1}+\sum_{i}u_{i}(x)h_{i}^{\vee},
n_{+}(x)p_{-n}n_{+}(x)^{-1}]\rangle ; 
$$
since by Lemma \ref{addlemma1bis},
$\pa_{n}(\sum_{i}u_{i}(x)h_{i}^{\vee})$ coincides with the Cartan part
of 
$$
-[p_{-1}, n_{+}(x)p_{-n}n_{+}(x)^{-1}],
$$
this is verified. 
So (\ref{evolphi}) holds, and the nonlocal parts of $\cH$ and $\cH_{0}$
coincide. 

So $\cH_{1}=\cH-\cH_{0}$ has no nonlocal terms, and satisfies
$\pa_{y}\cH_{1}=[\cH_{1}, p_{-1}]$; the arguments used to establish
(\ref{partialresult1}) then show that $\cH_{1}=0$. 
\hfill $\Box$

\medskip
\noindent
\begin{remark}
As we noted in the introduction, it is possible to derive
(\ref{PBn-n}) using Prop. \ref{prop3.1} and a differential equation in
the first variable, satisfied by $P(n_{+}\otimes n_{+})$. But we find
the present use of $B_{-}$ more natural.  \hfill $\Box$
\end{remark}

\subsection{Gelfand-Dickey-Dorfman structure}

According to Sect. 1.5, the VPA $\pi_0 = \CC[N_+/A_+]$ is endowed with a
Gelfand-Dickey-Dorfman structure. We give here a geometric interpretation
of it.

\begin{proposition} \label{propGDD}
Let $V((\la))$, $W((\mu))$ be $\G$-modules, and $v\in V$, $w\in W$
be such that ${\frak a_{+}}v=0$, ${\frak a_{+}}w=0$. We have:
$$
{\cal V}_{n_{+}v\otimes 1}(1\otimes n_{+}w)=
\sum_{k\ge 0}(-\pa)^{k}[(\Ad (n_{+})(t-a)_{k}^{(1)})_{+}n_{+}v]\otimes
(\Ad (n_{+})(t-a)_{k}^{(2)})_{+}n_{+}w
$$
\end{proposition}

{\em Proof.\/}
We apply the definition (\ref{hamvect}) of ${\cal V}_f$ to formula
(\ref{finalPBn-n}).
\hfill $\Box$

\subsection{Compatible nonlocal VPA structures}

Let us show how the nonlocal VPA structure on $\CC[B_{-}\times N_{+}]$,
defined in Thm. 3.1, can be embedded into an infinite family of
compatible nonlocal VPA structures (we call such a family compatible, if
any linear combination of these structures is again a nonlocal VPA
structure). 

Let us identify ${\frak g}$ with the subalgebra of the loop algebra
$\bar{\frak g}\otimes\CC((\la))$ of a
finite dimensional semisimple Lie algebra $\bar{\frak g}$, consisting of
the elements $x(\la)$ satisfying 
$$
x(\zeta\la)=x(\la)^{\sigma},
$$
$\sigma$ an automorphism of $\bar{\frak g}$ and $\zeta$ a root of
unity of the same order $r$. There is an action of $\CC((\la^{r}))$ on
${\frak g}$.  Let us denote by the same letter elements of
$\CC((\la^{r}))$ and the corresponding operators on ${\frak g}$.

For $n\in \ZZ$, set
$$
t_{n}=(\la^{rn}\otimes 1)t, \quad a_{n}=(\la^{rn}\otimes 1)a. 
$$
Note that $\ad p_{-1}$ commutes with the action of $\la^{rn}$, so that 
$$
t_{n}-a_{n}\in\Imm(\ad p_{-1})^{\otimes 2}. 
$$

Let us define an operation $P_n$ on $\CC[B_{-}\times N_{+}]$, by the
formula
\begin{equation}
\begin{array}{rcl}
P_{n}(g\otimes g)=[-\ell^{\otimes 2}(t_{n})(g(x)\otimes g(y))+r^{\otimes
2}(a_{n})(g(x)\otimes g(y))]\pa_{x}^{-1}\delta_{xy} \\ +
\sum_{n\ge 0}r^{\otimes 2}\Big((\ad p_{-1})^{-n-1}\otimes
1)(t_{n}-a_{n})\Big)(g(x)\otimes g(y))\pa_{x}^{n}\de_{xy}.
\label{compatPBg-g}
\end{array}
\end{equation}

\begin{proposition} \label{prop4.1}
The formulae (\ref{compatPBg-g}) define compatible nonlocal VPA
structures on $\CC[B_{-}\times N_{+}]$ (endowed with the derivation
$r(p_{-1})$). 
\end{proposition}

{\em Proof. \/} A combination of the brackets (\ref{compatPBg-g})
corresponds to the same formula, with $t_{n}$ and $a_{n}$ replaced by
$t_{f}$ and $a_{f}$ respectively, with $t_{f}=(f\otimes 1)t$ and
$a_{f}=(f\otimes 1)a$, $f$ a certain element of $\CC[\la^{r},\la^{-r}]$. 
The resulting bracket satisfies the conditions of Prop. \ref{prop1.1};
the elements of
$E_{i}$ in the second condition are $\ell^{\otimes 2}(t_{f})-r^{\otimes
2}(t_{f})$ for $i=0$,
and $0$ for $i>0$.  
\hfill $\Box$
\medskip

\noindent
\begin{remark} \label{rem4.1}
The variety $B_{-}\times N_{+}$ can be considered as an open subset of
$G$. It has a compatible family of nonlocal VPA structures defined by
(\ref{compatPBg-g}). In the same way as in the proof of
Prop. \ref{prop4.1}, 
one can show that formula (\ref{compatPBg-g}) defines a compatible
family of nonlocal VPA structures on the whole group $G$.
\hfill $\Box$
\end{remark}

\noindent
\begin{remark} \label{rem4.2}
The extension to $G$ of the nonlocal VPA structure defined in
Thm. \ref{thm3.1} is clearly left $G$-invariant. It follows that left
$G$-translations provide symmetries of the mKdV hierarchy, respecting
the Poisson structure. Infinitesimal left translations by elements of
${\frak n}_{+}$ correspond to the Toda flows; left translations by
elements of ${\frak b}_{-}$ do not change the variables $u_{i}$. A
class of translations that would be interesting to study further are
left translations by elements of the affine Weyl group; they should
mix local and nonlocal variables while respecting the Poisson
structure. A. Orlov pointed out to us that they probably coincide with
the Darboux transformations.
\hfill $\Box$
\end{remark}

\section{The proof of Lemma \ref{lemma3.5} and Lemma \ref{lemma3.6}}

\subsection{Proof of Lemma \ref{lemma3.5}} 
Equation (\ref{prePBb-n}) is rewritten as 
\begin{equation}
\begin{array}{rcl}
P(b_{-}\otimes n_{+})=-\sum_{\beta}\bar e_{\beta}b_{-}(x)
\otimes (\Ad (b_{-}(y)^{-1})
(\bar e^{\beta}))_{+}n_{+}(y)\pa_{x}^{-1}\de_{xy}
\\ +\sum_{n\in I}b_{-}(x)(\Ad (n_{+}(x))(p_{-n}))_{-}\otimes
n_{+}(y)p_{n}\pa_{x}^{-1}\delta_{xy}
+\on{\ local \ terms},
\end{array}
\label{PBb-n}
\end{equation}
with $\sum \bar e^{\beta}\otimes \bar e_{\beta}=\sum e^{\al}\otimes
e_{\al}+\sum_{i}h_{i}\otimes h^{\vee}_{i}$. The differential equations
satisfied by the left and right hand sides of (\ref{PBb-n}) are
$$
\begin{array}{rcl}
\pa_{y}(\on{lhs\ of\ (\ref{PBb-n})})=(\on{lhs\ of\ (\ref{PBb-n})})
(1\otimes p_{-1})
-1\otimes(p_{-1}+\sum_{i} u_{i}(y)h_{i}^{\vee})
\\
(\on{lhs\ of\ (\ref{PBb-n})})
-\sum_{i}\Ad(b_{-}(y))([p_{-1},h_{i}])b_{-}(x)
\otimes h_{i}^{\vee}n_{+}(y) \pa_{x}^{-1}\de_{xy} \\ +\on{\ local \
terms,}
\end{array}
$$
using (\ref{evoln+}) and (\ref{PBb-phi}),
and
$$
\begin{array}{rcl}
\pa_{y}(\on{rhs\ of\ (\ref{PBb-n})})=(\on{rhs\ of\ (\ref{PBb-n})})
(1\otimes p_{-1})
-1\otimes(p_{-1}+\sum_{i} u_{i}(y)h_{i}^{\vee}) \\ (\on{rhs\ of\ 
(\ref{PBb-n})})
  -\sum_{\beta}\bar
e_{\beta}b_{-}(x)\otimes \{[p_{-1}+  \sum_{i}h_{i}^{\vee} u_{i}(y), 
(\Ad (b_{-}^{-1}(y))(\bar
e^{\beta}))_{+}] \\ -([p_{-1}+\sum_{i}h_{i}^{\vee} u_{i}(y), \Ad
(b_{-}^{-1}(y))
(\bar e^{\beta})])_{+}\} \\ n_{+}(y)\pa_{x}^{-1}\de_{xy}+\on{\ local \
terms,}
\end{array}
$$
using (\ref{evoln+}) and (\ref{evolb-}). 
Since 
\begin{equation}
\sum_{i}\Ad(b_{-}(y))([p_{-1},h_{i}])\otimes h_{i}^{\vee}=\sum_{\beta}
\bar
e_{\beta}\otimes [p_{-1},(\Ad (b_{-}(y)^{-1})(\bar e^{\beta}))_{1}],
\label{aux1}
\end{equation}
these two equations coincide. In (\ref{aux1}), we denote by $x_{1}$ the
part of $x\in{\frak g}$, of principal degree 1. (\ref{aux1}) is proved by
pairing its right and left hand sides with $1\otimes h_{i}$,
$i=1,\ldots,l$. It follows that the difference of the two sides of
(\ref{PBb-n}) satisfies 
\begin{equation}
\begin{array}{rcl}
\pa_{y}(\on{ lhs\ of\ (\ref{PBb-n})}-\on{ rhs\ of\ (\ref{PBb-n})})
=(\on{ lhs\ of\ (\ref{PBb-n})}  
-\on{ rhs\ of\ (\ref{PBb-n})})(1\otimes p_{-1})
\\
-(1\otimes 
(p_{-1}+\sum_{i} u_{i}(y)h_{i}^{\vee}))(\on{ lhs\ of\ (\ref{PBb-n})}
-\on{ rhs\ of\ (\ref{PBb-n})})+\on{local\ terms.}
\label{evidence1} 
\end{array}
\end{equation}

On the other hand, we have 
$$
\begin{array}{rcl}\pa_{x}(\on{ lhs\ of\ (\ref{PBb-n})})=(\on{ lhs\ of\ 
(\ref{PBb-n})})((p_{-1}+\sum_{i}
 u_{i}(x)h_{i}^{\vee})\otimes 1) \\ -
\sum_{i,k}[b_{-}(x)h_{i}^{\vee}\pa_{k} u_{i}(x)\otimes n_{+}(y)p_{k}]
\pa_{x}^{-1}\de_{xy}
 +\on{ \ local\ terms,}
\end{array}
$$
and
$$
\begin{array}{rcl}\pa_{x}(\on{ rhs\ of\ (\ref{PBb-n})})=(\on{ rhs\ of\ 
(\ref{PBb-n})})((p_{-1}+\sum_{i}
 u_{i}(x)h_{i}^{\vee})\otimes 1)
 \\ 
+\sum_{n\in I}b_{-}(x) \{[p_{-1}+\sum_{i}h_{i}^{\vee} u_{i}(x), (\Ad
(n_{+}(x))(p_{-n}))_{-}]
\\
-([p_{-1}+\sum_{i}h_{i}^{\vee} u_{i}(x),\Ad
(n_{+}(x))(p_{-n})])_{-} \\ \otimes n_{+}(y)p_{n}\pa_{x}^{-1}\delta_{xy}
\}
+\on{ \ local\ terms.}
\end{array}
$$

We have the equality
$$
[p_{-1},(\Ad (n_{+}(x))(p_{-n}))_{1}]=\sum_{i}h_{i}^{\vee}\pa_{n}
u_{i}(x), 
$$
because of Lemma \ref{addlemma1bis}. 
 
Therefore the right hand sides of the last two formulas coincide up to
local terms, and 
\begin{equation}
\begin{array}{rcl}
\pa_{x}(\on{ lhs\ of\ (\ref{PBb-n})}-\on{ rhs\ of\ (\ref{PBb-n})})
=(\on{ lhs\ of\ (\ref{PBb-n})}-\on{ rhs\ of\ (\ref{PBb-n})}) \\
(1\otimes 
(p_{-1}+\sum_{i} u_{i}(x)h_{i}^{\vee}))+\on{local\ terms}.
\label{evidence2} 
\end{array}
\end{equation}

Let $\cX=(1\otimes n_{+}(y)^{-1})(\on{lhs\ of\ (\ref{PBb-n})}-\on{rhs\
of\ (\ref{PBb-n})})(1\otimes b_{-}(x)^{-1})$, then 
$$
\pa_{x}\cX=\on{local\ terms}
$$
by (\ref{evidence2}), and 
$$
\pa_{y}\cX=[\cX, 1\otimes p_{-1}]+\on{local\ terms}
$$
by (\ref{evoln+}) and (\ref{evidence1}). The first equation gives 
$$
\pa_{x}\cX=\sum_{n\ge 0}\cX_{n}(y)\pa_{x}^{n}\de_{xy},
$$
so 
$$
\cX=\cX_{0}(y)\pa_{x}^{-1}\de_{xy}+\on{local\ terms,}
$$
and the second equation gives us
$$
\pa_{y}\cX_{0}(y)=[\cX_{0}(y), 1\otimes p_{-1}].
$$
Let $\xi$ be any element of the dual to ${\frak b}_{-}$, and
$\cX_{\xi}(y)=(\xi\otimes 1)(\cX_{0}(y))$. Then $\cX_{\xi}(y)$ has
values in ${\frak n}_{+}$ and satisfies
\begin{equation}
\pa_{y}\cX_{\xi}(y)=[\cX_{\xi}(y),p_{-1}]; 
\label{evolX}
\end{equation}
we then follow the proof of \cite{EFr}, lemma 3, to conclude that 
$\cX_{\xi}(y)$ is
constant and lies in ${\frak a}_{+}$. Recall how this can be done: 
decompose $\cX_{\xi}(y)$ in its homogeneous principal components
$\sum_{i}\cX_{\xi,i}(y)$, and each component along the decomposition
$\Imm(\ad p_{-1})\oplus{\frak a}$, as
$\cX_{\xi,i}^{1}(y)+\cX_{\xi,i}^{2}(y)$; let $i$ be the smallest index,
such that $\cX_{\xi,i}^{1}(y)$ is not zero; the equation implies that
$\pa_{y}\cX_{\xi}(y)$ has a nonzero component of degree $i-1$ in
$\Imm\ad p_{-1}$, hence a contradiction. So $\cX_{\xi}(y)$ lies in
${\frak a}$; (\ref{evolX}) then implies that it is constant. We finally 
obtain: 
$$
\on{ lhs\ of\ (\ref{PBb-n})}-\on{ rhs\ of\ (\ref{PBb-n})}
=\sum_{n\in I}
x_{n}b_{-}(x)\otimes
n_{+}(y)p_{n}\pa_{x}^{-1}\de_{xy}+\on{\  local\ terms,}
$$
with $x_{n}\in{\frak b}_{-}$. But there is only one possibility,
$x_{n}=0$, which is compatible with the following invariance property
of $P$.

Recall that for $\xi\in{\frak b}_{-}$, $\ell(\xi)$ is the derivation
of the algebra $\bar\pi_{0}$ defined by the action of the left
translation by $\xi$, on $B_{-}\times N_{+}$. Since $\ell(\xi)$
commutes with $\pa$, it induces an endomorphism (also denoted by
$\ell(\xi)$) of $\cR_{2}(\bar\pi_{0})$, according to the rules used in
2.1 in the case of $\pa_{n}$. We then have:

\begin{lemma} \label{lemma3.2} 
For $a,b\in\bar\pi_{0}$, $\xi\in{\frak b}_{-}$,
\begin{equation}
P(\ell(\xi)a\otimes b)+P(a\otimes \ell(\xi)b)=\ell(\xi)P(a\otimes b).
\label{invariance}
\end{equation}
\end{lemma}
{\em Proof\/}. For the brackets $P(b_{-}\otimes b_{-})$, this follows
from (\ref{PBb-b}) and the invariance of $t$. We also have
$$
\begin{array}{rcl}
P(\ell(\xi)b_{-}\otimes u_{i})=
\pa_{y}([\xi,\Ad(b_{-}(y))(h_{i})]b_{-}(x)\pa_{x}^{-1}\de_{xy}
  +\Ad(b_{-}(y))(h_{i})\xi b_{-}(x)\\ \pa_{x}^{-1}\de_{xy})=\ell(\xi)
P(b_{-}\otimes u_{i})
\end{array}
$$
so that 
$$
P(H_{n}\otimes\ell(\xi)b_{-})=\ell(\xi)P(H_{n}\otimes b_{-}),
$$
(in this equality, the second $\ell(\xi)$ is $\ell(\xi)\otimes 1$ acting
on $\End(V)((\xi))\otimes \cR_{2}(\bar\pi_{0})$)
and $(\ell(\xi)\otimes 1)r_{n}=(1\otimes\ell(\xi))r_{n}$ (equality in
$\End(V)((\xi))\otimes \cR_{2}(\bar\pi_{0})$; so that 
$$
P(F_{n}\otimes\ell(\xi)b_{-})=\ell(\xi)P(F_{n}\otimes b_{-}).
$$

Finally, the elements of $\pi_{0}^{+}$ are invariant under $\ell(\xi)$, so
the identity is trivially satisfied for their Poisson brackets. 
\hfill $\Box$
\medskip

Now (\ref{PBb-n}) follows.
\hfill $\Box$

\subsection{Proof of Lemma \ref{lemma3.6}}

The operation $P_{0}$ is defined by the identities
$$
P_{0}(b_{-}\otimes b_{-})=-\ell^{\otimes 2}(t)(b_{-}(x)\otimes b_{-}(y))
\pa_{x}^{-1}\de_{xy}, 
$$
\begin{equation}
\begin{array}{rcl}
P_{0}(b_{-}\otimes n_{+})=[-b_{-}(b_{-}^{-1}t^{(1)}b_{-})_{-}(x)
\otimes (b_{-}^{-1}t^{(2)}b_{-})_{+}(y)n_{+}(y)
\\
+b_{-}(n_{+}a^{(1)}n_{+}^{-1})_{-}(x)\otimes
(n_{+}a^{(2)}n_{+}^{-1})_{+}n_{+}(y)] \pa_{x}^{-1}\de_{xy}
\\
+\sum_{n\ge 0}[b_{-}(n_{+}(t-a)_{n}^{(1)}n_{+}^{-1})_{-}(x)\otimes 
(n_{+}(t-a)^{(2)}_{n}n_{+}^{-1})_{+}n_{+}(y)]\pa^{n}\de_{xy},
\label{attemptPBb-n}
\end{array}
\end{equation}
where we denote $((\ad\ p_{-1})^{-n-1}\otimes 1)(t-a)$ by $(t-a)_{n}$,
any element $\al\in{\frak g}\otimes {\frak g}$ is decomposed as
$\sum \al^{(1)}\otimes \al^{(2)}$, and 
\begin{equation}
\begin{array}{rcl}
P_{0}(n_{+}\otimes n_{+})=r^{\otimes
2}(a)(n_{+}(x)\otimes n_{+}(y))\pa_{x}^{-1}\delta_{xy} \\ +
\sum_{n\ge 0}r^{\otimes 2}\Big((\ad p_{-1})^{-n-1}\otimes
1)(t-a)\Big)(n_{+}(x)\otimes n_{+}(y))\pa_{x}^{n}\de_{xy}.
\end{array}
\label{attemptPBn-n}
\end{equation}
By construction, $P_{0}(g\otimes g)$ satisfies the identities 
\begin{equation}
\pa_{x}P_{0}(g\otimes g)=(r(p_{-1})\otimes 1)P_{0}(g\otimes g), \quad
\pa_{y}P_{0}(g\otimes g)=(1\otimes r(p_{-1}))P_{0}(g\otimes g). 
\label{evolPzero}
\end{equation}

Clearly, $P_{0}(b_{-}\otimes b_{-})$ coincides with $P(b_{-}\otimes
b_{-})$. $\cB=P(b_{-}\otimes n_{+})$ satisfies the equation
$$
\begin{array}{rcl}
\pa_{y}\cB+(1\otimes(p_{-1}+\sum_{i}u_{i}h_{i}^{\vee}))\cB-\cB(1\otimes
p_{-1})=
\sum_{i} 
\left( \pa_{y} \Ad (b_{-}(y))(h_{i})b_{-}(x)\pa_{x}^{-1} \right. \\
\otimes
\left.  h_{i}^{\vee}n_{+}(y) \right) \de_{xy},
\end{array}
$$
by virtue of (\ref{PBb-phi}) and (\ref{evoln+}). Let us determine an
equation satisfied by $\cB_{0}=P_{0}(b_{-}\otimes n_{+})$. Let
$\cE=P_{0}(b_{-}\otimes b_{-})$. $\cE$ satisfies the equation 
$$
\pa_{y}\cE=\cE(1\otimes (p_{-1}+\sum_{i}u_{i}h_{i}^{\vee}))-\pa_{y}[\Ad
(b_{-}(y))(h_{i})b_{-}(x)\pa_{x}^{-1} \otimes b_{-}(y)] \de_{xy}.
$$
Note that due to (\ref{evolPzero}), $P_{0}(b_{-}\otimes g)$ satisfies 
$$
\pa_{y}P_{0}(b_{-}\otimes g)=(1\otimes r(p_{-1}))P_{0}(b_{-}\otimes g). 
$$
This implies, writing $\cB_{0}$ as $P_{0}(b_{-}\otimes
b_{-}^{-1}g)$ [$B_{-}\times N_{+}$ has a left $B_{-}$-action, defined as
the product of the left action of $B_{-}$ on itself and of the trivial
one, that we use here], that $\cB_{0}$ satisfies 
$$
\begin{array}{rcl}
\pa_{y}\cB_{0}+(1\otimes(p_{-1}+\sum_{i}u_{i}h_{i}^{\vee}))\cB_{0}
-\cB_{0}(1\otimes
p_{-1})=
\sum_{i} \left(
\pa_{y} \Ad (b_{-}(y))(h_{i})b_{-}(x)\pa_{x}^{-1} \right. \\
\left. \otimes 
h_{i}^{\vee}n_{+}(y) \right) \de_{xy}.
\end{array}
$$
Let us set $\cB_{1}=(1\otimes n_{+}^{-1}(y))(\cB-\cB_{0})$, we obtain 
$$
\pa_{y}\cB_{1}+[1\otimes p_{-1}, \cB_{1}]=0.
$$
The nonlocal parts of $\cB$ and $\cB_{0}$ coincide, so that $\cB_{1}$
contains only local terms; write $\cB_{1}=\sum_{n\ge
0}\cB_{1}^{(n)}(x)\pa_{x}^{n}\de_{xy}$ (each $\cB_{1}^{(n)}$ belongs to
the tensor product of the tangent space $T_{b_{-}(x)}B_{-}$
to $B_{-}$ at $b_{-}(x)$ with 
${\frak n}_{+}$), we then get 
$$
[1\otimes p_{-1}, \cB_{1}^{(0)}(x)]=0, \quad
\cB_{1}^{(n)}=[1\otimes p_{-1}, \cB_{1}^{(n+1)}(x)]
$$
for $n\ge 0$, so $\cB_{1}^{(0)}(x)$ belongs to $T_{b_{-}(x)}B_{-}\otimes
{\frak
a}$ by the first equation and to $T_{b_{-}(x)}B_{-}\otimes \Imm (\ad\
p_{-1})$ by the second one (specialized to $n=0$), so that it is zero; 
repeating
this argument for $\cB_{1}^{(1)}$, we find it to vanish as well, etc. So
$\cB_{1}=0$ and 
\begin{equation}
P(b_{-}\otimes n_{+})=P_{0}(b_{-}\otimes n_{+}).
\label{partialresult1}
\end{equation}

Now, $\cB=P(b_{-}\otimes n_{+})$ satisfies the equation
\begin{equation}
\pa_{x}\cB=\cB((p_{-1}+\sum_{i}u_{i}(x)h_{i}^{\vee})\otimes
1)+\sum_{i}b_{-}(x)\otimes (h_{i}^{\vee}P(u_{i}\otimes n_{+})).
\label{evolBinx}
\end{equation}

On the other hand, let $\cC=P(n_{+}\otimes n_{+})$ and
$\cC_{0}=P_{0}(n_{+}\otimes n_{+})$. $\cC$ satisfies the equation 
$$
\pa_{x}\cC=-((p_{-1}+\sum_{i}u_{i}(x)h_{i}^{\vee})\otimes
1)\cC+\cC(p_{-1}\otimes 1)-\sum_{i}h_{i}^{\vee}n_{+}(x)\otimes
P(u_{i}\otimes n_{+}).
$$
Let us determine an equation satisfied by $\cC_{0}$. Due to
(\ref{evolPzero}), $P_{0}(g\otimes n_{+})$ satisfies 
$$
\pa_{x}P_{0}(g\otimes n_{+})=P_{0}(g\otimes n_{+})(p_{-1}\otimes 1), 
$$
and writing $P_{0}(n_{+}\otimes n_{+})$ as $P_{0}(b_{-}^{-1}g\otimes
n_{+})$ (using the same left $B_{-}$-action as above) and using
(\ref{evolBinx}), we get 
$$
\pa_{x}\cC_{0}=-((p_{-1}+\sum_{i}u_{i}(x)h_{i}^{\vee})\otimes
1)\cC_{0}+\cC_{0}(p_{-1}\otimes 1)-\sum_{i}h_{i}^{\vee}n_{+}(x)\otimes
P(u_{i}\otimes n_{+}).
$$
$\cC$ and $\cC_{0}$ satisfy the same equation, so that
$\cC_{1}=(n_{+}(x)^{-1}\otimes 1)(\cC-\cC_{0})$
(which belongs to the
tensor product of ${\frak n}_{+}$ with the tangent space to $N_{+}$ at
$n_{+}(y)$) 
satisfies 
$$
\pa_{x}\cC_{1}=[\cC_{1}, p_{-1}\otimes 1]. 
$$
Since the nonlocal parts of $P_{0}$ and $P$ coincide, $\cC_{1}$ contains
no nonlocal terms. We can use the same arguments as in the case of
$\cB_{1}$, to conclude that $\cC_{1}=0$ and 
\begin{equation}
P(n_{+}\otimes n_{+})=P_{0}(n_{+}\otimes n_{+}). 
\label{partialresult2}
\end{equation}

Lemma \ref{lemma3.6} now follows from (\ref{partialresult1}),
(\ref{partialresult2}).
\hfill $\Box$

\end{document}